\documentclass[english]{article}
\usepackage[T1]{fontenc}
\usepackage[latin9]{inputenc}
\usepackage{array}
\usepackage{url}
\usepackage{graphicx}
\usepackage{setspace}
\usepackage{amssymb}

\makeatletter

\providecommand{\tabularnewline}{\\}

\newenvironment{lyxlist}[1]
{\begin{list}{}
{\settowidth{\labelwidth}{#1}
 \setlength{\leftmargin}{\labelwidth}
 \addtolength{\leftmargin}{\labelsep}
 }}
{\end{list}}

\makeatother

\usepackage{babel}

\begin{document}

\title{\textbf{10 conjectures in additive number theory}}

\author{Benoit Cloitre }
\maketitle
\begin{abstract}
Following an idea of Rowland \cite{key-1} we give a conjectural way
to generate increasing sequences of primes using {}``gcd-algorithms''.
These algorithms seem not so useless for searching primes since it
appears we found sometime primes much more greater than the number
of required iterations. In an other hand we propose new formulations
of famous conjectures from the additive theory of numbers (the weak
twin prime conjecture, the Polignac conjecture, the Goldbach conjecture
or the very general Schinzel's hypothesis H). For the moment these
are experimental results obtained using pari-gp \cite{key-4}. 
\end{abstract}

\section*{Introduction}

In \cite{key-1} the author solved a question related to the sequence
defined recursively by $R(1)=7$ and for $n\geq2$ by:
\begin{itemize}
\item $R(n)=R(n-1)+\gcd(R(n-1),n)$ 
\end{itemize}
Namely he showed that $R(n)-R(n-1)$ is always a prime number or 1.
At first glance this recursion is not very useful for finding primes
since these differences generate primes quite randomly and the prime
values are less than the number of required iterations. However, looking
closer to Rowland recursion we found a way to exhibit the records
values in the sequence of differences $R(n)-R(n-1)$ (these differences
are given by $A137613$ in \cite{key-2}). Namely $n+1$ when $R(n)=2n+2$
are exaclty those records values. So this gives a method to exhibit
primes in increasing order (and the growth is exponential). Those
$n+1$ such that $R(n)=2n+2$ begin:
\begin{lyxlist}{00.00.0000}
\item [{$5,11,23,47,101,233,467,941,1889,3779,7559,15131,30323,60647,\ldots$$ $}]~
\end{lyxlist}
Although this seems provable using Rowland reasoning, it looks like
it is unknown since it is not mentioned before in \cite{key-1} and
this sequence of records values is not in \cite{key-2}. This is clearly
a nice fact and we try to generalise this observation. In order to
facilitate the computations we use the absolute value and run the
algorithm {}``backwards'' as we will see. So we succeeded to extend
somewhat this result with the conjectures 1 and 2 and propose a method
for building primes (section 1). In a conjecture 3 (section 2) we
claim there are sequences of lesser of twin primes growing very fast.
Then we provide a conjecture 3bis (cf. 2.4) related to the Polignac
conjecture and a conjecture 3ter (cf. 2.5) dealing with the prime
triplet conjecture. In a conjecture 4 (section 3) we propose a way
to find very big primes compared to the number of iterations which
could improve the conjecture 2 and be an efficient tool for searching
primes. In a conjecture 5 (section 4) we propose a first constructive
way to prove the Goldbach conjecture but with some restiction due
to an exceptional set of non working cases. In a conjecture 6 (section
5) we propose to summarize all our observations with a reformulation
of the Schinzel's hypothesis H. The conjecture 7 (cf. 5.2) gives another
approach of the conjecture 6.

We also discuss the Shevelev conjecture who extends Rowland idea for
twin primes. Indeed, V. Shevelev \cite{key-3} introduced the related
sequence ($A166944$ in \cite{key-2}) :
\begin{itemize}
\item $S(1)=2$ 
\item $S(n)=S(n-1)+\gcd(S(n-1),n-1+(-1)^{n})$
\end{itemize}
And noticed that the records values of differences ($A166945$ in
\cite{key-2}) :
\begin{lyxlist}{00.00.0000}
\item [{$7,13,43,139,313,661,1321,...$}]~
\end{lyxlist}
gives for terms $>7$ the greater prime of some twin primes pairs
. As for Rowland sequence it can be seen that the indices where these
records occur is given by the sequence of $n+1$ such that $S(n)=2n+1$
(except for $7$) and $(n,S(n)-S(n-1))$ is then a twin pair of primes.
Using this algorithm backwards the conjecture 8 (section 6) propose
a sligthly different way than the conjecture 3 or 6 or Shevelev conjectures
to prove the weak twin prime conjecture. 

Although Shevelev managed to do it, it is not easy to generalise this
observation using the {}``forward'' original Rowland recursion since
we need to find where the records occur. However using the absolute
value we shall see it is easy to obtain many increasing sequences
of twin primes or of triplet of primes since we just have to check
indices where zeros occur. Thus in a conjecture 9 (section 7) we propose
a method for generating primes $m$-uplet of any type using periodic
sequences in the gcd algorithm. Then we extend the idea to a family
of polynomials giving another version of the Schinzel's hypothesis
H. This allows us to perform easier computations than using the conjecture
6 for searching $m$-uplet.

We finally merge our {}``backward'' recursion and Shevelev idea
for approaching the Goldbach conjecture in the conjecture 10 (section
8) in a nicer way than the conjecture 5. There is apparently no more
exception for $n$ large enough. This right way to deal with this
old conjecture using gcd-recursion was very hard to find. We also
discuss a less known but hard conjecture of Legendre (section 9).

There is a curious fact in this study. It appears our methods work
often better for large integer since once we found a good starting
value we generate only primes, despite the probability to pick up
big primes among large integers goes to zero.

\section{Variation on Rowland recursion}

After some attempts we arrive to this recursion seeming generating
primes in a general way starting always with the same initial value
$1$. For a given integer value $m\geq1$ we define the sequence $\left(a(n)\right)_{n\geq1}$
by $a(1)=1$ and for $n\geq2$ by the recursion:
\begin{itemize}
\item $a(n)=\left|a(n-1)-\gcd\left(a(n-1),mn-1\right)\right|$.
\end{itemize}
Now let us consider the values of $n$ such that we get:
\begin{itemize}
\item $a(n)=0\Leftrightarrow$ $a(n+1)=mn+m-1$. 
\end{itemize}
We claim that $\forall m\geq1$ this sequence of indices $n$ gives
rise to an infinite sequence $\left(b_{m}(k)\right)_{k\geq1}$.

\subsubsection*{Example}

For $m=1$ the sequence$\left(a(n)\right)_{n\geq1}$ begins:
\begin{lyxlist}{00.00.0000}
\item [{$1,0,2,1,0,5,4,3,2,1,0,11,10,9,8,7,6,5,4,3,2,1,0,23,22,21,20,19,\ldots$}]~
\end{lyxlist}
And the sequence $\left(b_{1}(k)\right)_{k\geq1}$ of indices where
zeros appear in $\left(a(n)\right)_{n\geq1}$ begins:
\begin{lyxlist}{00.00.0000}
\item [{$2,5,11,23,47,79,157,313,619,1237,2473,4909,9817,19603,39199,78193,\ldots$}]~
\end{lyxlist}
$ $These listed numbers are prime numbers. This led us to a first
conjecture.

\subsection{Conjecture 1}

We claim:
\begin{itemize}
\item $b_{1}(k)$ is prime for $k\geq1$ and $b_{1}(k)\sim c2^{k}$ $(k\rightarrow\infty)$
with $c=1.186\ldots$
\end{itemize}
See the APPENDIX 1 for a table supporting this conjecture. Although
this conjecture is similar to the record sequence mentioned in the
introduction, it is more interesting to us. Indeed we were able to
generalise the result and much more seems true. So we make a stronger
conjecture.

\subsection{Conjecture 2}

In general we claim that for any $m\geq1$:
\begin{itemize}
\item $mb_{m}(k)+m-1$ is prime for $k$ large enough (usually $k\geq2$
is working for small $m$). 
\item $ $$mb_{m}(k)+m-1\sim c_{m}(m+1)^{k}$ $(k\rightarrow\infty)$ with
$c_{m}>0$. 
\end{itemize}

\subsubsection*{Additional examples supporting the conjecture 2}

For $m=3$ the sequence $\left(a(n)\right)_{n\geq1}$ begins: 
\begin{lyxlist}{00.00.0000}
\item [{$1,0,8,7,0,17,16,15,14,13,12,11,10,9,8,7,6,5,4,3,2,1,0,71,\ldots$}]~
\end{lyxlist}
And the sequence $\left(b_{3}(k)\right)_{k\geq1}$ begins: 
\begin{lyxlist}{00.00.0000}
\item [{$2,5,23,89,337,1335,5307,\ldots$}]~
\end{lyxlist}
Next the values $3b_{3}(k)+2$ appear to be prime values for $k\geq2$: 
\begin{lyxlist}{00.00.0000}
\item [{$8,17,71,269,1013,4007,15923,63521,253949,1014317,\ldots$}]~
\end{lyxlist}
And $3b_{3}(n)+2\sim c_{3}4^{n}$ $(n\rightarrow\infty)$ with $c_{3}=0.96\ldots$

\begin{flushleft}
For $m=4$ the sequence $\left(a(n)\right)_{n\geq1}$ begins:
\par\end{flushleft}
\begin{lyxlist}{00.00.0000}
\item [{$1,0,11,10,9,8,7,6,5,4,3,2,1,0,59,58,57,\ldots$}]~
\end{lyxlist}
And the sequence $\left(b_{4}(k)\right)_{k\geq1}$ begins:
\begin{lyxlist}{00.00.0000}
\item [{$2,14,62,314,1574,7846,38020,\ldots$}]~
\end{lyxlist}
Next the values $4b_{4}(k)+3$ appear to be prime values only for
$k\geq1$:
\begin{lyxlist}{00.00.0000}
\item [{$11,59,251,1259,6299,31387,152083,758971,3790651,18953251,\ldots$}]~
\end{lyxlist}
And $4b_{4}(n)+3\sim c_{4}5^{n}$ $(n\rightarrow\infty)$ with $c_{4}=1.9408\ldots$

\subsection{An efficient algorithm for finding primes?}

Note we find sometime primes greater than $n$ after making $n$ iterations.
So it could be an efficient method for finding primes since the gcd
algorithm is well known and {}``fast'' computation can be performed.
For instance if $m=28$ we compute $1000000$ terms of the sequence
$\left(a(n)\right)_{n\geq1}$ . In this range $a(n)$ vanishes $5$
times. This allows us to compute $\left(b_{28}(n)\right)_{1\leq n\leq5}$
wich gives $5$ values of $28b_{28}(n)+27$:
\begin{lyxlist}{00.00.0000}
\item [{$83,1147,31891,924811,26819491.$}]~
\end{lyxlist}
All these values are primes and the last one $26819491$ gives a prime
number larger than the $1000000$ iterations (see APPENDIX 2 for more
exemples). Moreover it appears this kind of algorithm can be adapted
for finding bigger primes as shown thereafter.

\subsubsection{A simple rule of construction}

Observe from our definition we have for $n\geq2$ (letting $b=b_{m}$
):
\begin{itemize}
\item $a(b(n)+1)=(m+1)a(b(n-1)+1)+m+m\sum_{j=b(n-1)+1}^{b(n)-1}\left(a(j+1)-a(j)+1\right)$ 
\end{itemize}
Hence we have an {}``almost'' recurrence relation between two consecutive
records values which are conjectured to be prime for $n$ large enough.
This formula explains also why these record values are growing like
$(m+1)^{n}$ since it appears we usually don't need to compute all
terms in the sum. Indeed experiments show that $a(j+1)-a(j)=-1$ for
$b(n-1)<n_{0}\leq j\leq b(n)-1$ and this $n_{0}$ depends on $b(n-1)$
and stays {}``near'' from this value. i.e. we claim that $n_{0}-b(n-1)\ll\sqrt{b(n-1)}$
and is often much more smaller. So we can launch a computation and
when the computer returns $a(j+1)-a(j)=-1$ sufficiently {}``often''
we may stop the computation and suspect we are beyond this $n_{0}$.
Therefore we built perhaps a prime greater than a starting prime.
It could then be interesting to know rules in order to choose {}``good''
values of $m$ making $n_{0}$ the closest of $b(n-1)$ as possible.

\subsubsection{Exemples}

For instance consider $m=10$ and the sequence of $a_{n}-a_{n-1}$
when $\left|a_{n}-a_{n-1}\right|>1$ wich yields somewhere for $a_{n}-a_{n-1}$:
\begin{itemize}
\item $43213789,-3,-13,-15241,-43,-1889,-3,-433,-113,-3,-5827,-247$
\end{itemize}
\begin{singlespace}
The positive terms are our record values $a(b(k)+1)$ and are primes
values. Here the computation returns nothing more using $10^{5}$
iterations. Hence we are perhaps beyond our $n_{0}$ and the next
record would be given by $11\times43213789+10-10\times(2+12+15240+42+1888+2+432+112+2+5826+246)=475113649$
which is a prime number and, as expected, our next record value. So
we have found a prime number 11 times greater than the given prime
$43213789$ with few iterations. 
\end{singlespace}

But let us see this with a more striking exemple and take $m=100000$.
We compute $\left|a_{n}-a_{n-1}\right|>1$ for $n\leq10^{6}$ wich
yields for these $a_{n}-a_{n-1}$:
\begin{itemize}
\item $299999,-59,29994499999,-3,-7,-53,-3$
\end{itemize}
So we suspect we have nothing more after $-3$ until the next record
value. Thus $100001\times29994499999+100000-100000\times(2+6+52+2)=2999479988299999$
should be the next record and it is indeed a prime value. So with
few effort we found a prime number $100001$ times greater than another
one. A simple routine under pari-gp found easily big primes with this
method (see the end of the APPENDIX 2).

\subsubsection{Probability to got a prime value from the first record}

We also observe something which could have practical use for searching
primes if people work together like for the GIMPS where they share
computer power. We keep previous notations and our definition of $a_{n}$
(depending on $m$) so that the first record is simply $3m-1$. Then
we tried to estimate the chance to got a prime value with the second
record value for various $m$ after making $\left\lfloor m^{\alpha}\right\rfloor $
iterations and where $0<\alpha\leq1$ . i.e. for a given $m$ we define:
\begin{itemize}
\item $f_{\alpha}(m)=1$ if the second record value which equals $(m+1)(3m-1)+m+m\sum_{j=4}^{\left\lfloor m^{\alpha}\right\rfloor }\left(a(j+1)-a(j)+1\right)$
is a prime value and otherwise $f_{\alpha}(m)=0$.
\item $L(\alpha)=\lim_{n\rightarrow\infty}\frac{1}{n}\sum_{k=1}^{n}f_{\alpha}(k)$ 
\end{itemize}
Then we conjecture that:
\begin{itemize}
\item $L(\alpha)$ exists with simply $L(\alpha)=\alpha$ and so we know
the probability to got a prime with this method.
\end{itemize}
To see this last point with experiments suppose $E$ is a set of $30000$
random values $m$ satisfying $10^{14}<m<10^{15}$ . We take $\alpha=\frac{1}{7}$
so that we have to make around $100$ iterations only to got the value
$(m+1)(3m-1)+m+m\sum_{j=4}^{\left\lfloor m^{\alpha}\right\rfloor }\left(a(j+1)-a(j)+1\right)$.
We then launch the computation 3 times (so that E changes each time)
and we plot the different graph of $\frac{1}{n}\sum_{i=1}^{n}f_{\alpha}(m_{i})$
(graphs are black, blue and green) where $m_{1}<m_{2}<...\in E$ compared
to the graph of $y=\frac{1}{7}$ (red).

\begin{center}
(fig.1)
\par\end{center}

\begin{center}
\includegraphics[width=0.5\paperwidth]{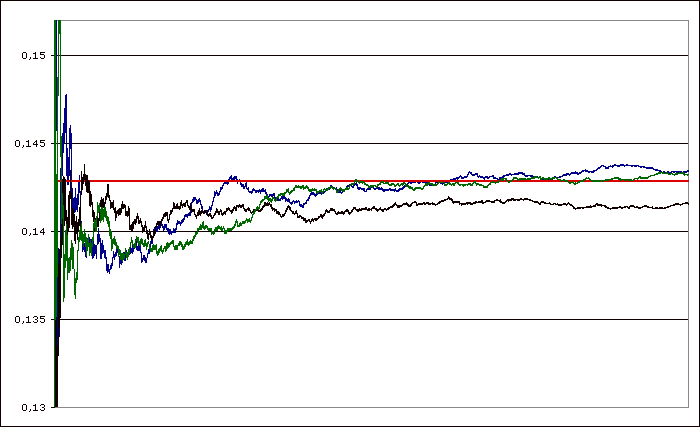}
\par\end{center}

This supports the claim $L(\alpha)=\alpha$. In other word this means
that if we take $\left\lfloor \alpha^{-1}\right\rfloor $ random values
of $m$ we are fairly certain the formula $(m+1)(3m-1)+m+m\sum_{j=4}^{\left\lfloor m^{\alpha}\right\rfloor }\left(a(j+1)-a(j)+1\right)$
will produce at least a prime number among them. Note we may start
form the second record value (or the third...) instead of the first
but we aim to work with big values of $m$ so it seems inappropriate
to compute the second record.

\subsection*{Remark}

Although this method for finding primes seems satisfying it is not
evident to see whether it is very efficient since we need to combine
two computations (one for finding the candidates and one for testing
primalty) and so far we find no rule in order to force $n_{0}$ to
stay very close from $b(n-1).$ Perhaps the conjecture 4 would be
better for that purpose (see 3.4.). By the way this kind of algorithms
is interesting in its own since variations of the Rowland-Shevelev
algorithm led us to formulate in a new way old conjectures from the
additive theory of numbers. Hereafter we come across conjectures like
the Polignac conjecture, the Schintzel hypothesis H or the Goldbach
conjecture and our approach suggests analytic or probabilistic study.

\section{An increasing sequence of twin primes}

This was somewhat surprising to find such a sequence since it is not
known whether there are infinitely many twin primes. Here we consider
the recursion $a(1)=1$ and for $n\geq2$:
\begin{itemize}
\item $a(n)=\left|a(n-1)-\gcd(a(n-1),n^{p}-1)\right|$
\end{itemize}

\subsection{Conjecture 3 }

Suppose $p\geq2$ is a prime number then we claim:
\begin{itemize}
\item there are infinitely many values of $n$ such that $a(n)=0$.
\item for $n$ large enough (usually $n>2$ is working for small $p$) we
have $a(n)=0\Rightarrow n$ is prime and $\frac{(n+1)^{p}-1}{n}$
is prime.
\end{itemize}

\subsection{Corollary}

There are infinitely many twin primes since the conjecture yields
for $p=2$ :
\begin{itemize}
\item for $n$ large enough we have $a(n)=0\Rightarrow n$ is prime and
$n+2$ is prime 
\end{itemize}
and this happens infinitely many times.

\subsection{Tables \textmd{\normalsize (in the sequel $\delta(n)=1$ if $n$
is prime and $0$ otherwise.)}}

\begin{center}
$p=2$
\par\end{center}

\begin{center}
\begin{tabular}{|c|c|c|c|}
\hline 
Values of $n$ such that $a(n)=0$ & $n+2$ & $\delta(n)$ & $\delta(n+2)$\tabularnewline
\hline
\hline 
2 & 4 & 1 & 0\tabularnewline
\hline 
11 & 13 & 1 & 1\tabularnewline
\hline 
137 & 139 & 1 & 1\tabularnewline
\hline 
19181 & 19183 & 1 & 1\tabularnewline
\hline 
367953497 & 367953499 & 1 & 1\tabularnewline
\hline
\end{tabular}
\par\end{center}

\begin{center}
$p=3$
\par\end{center}

\begin{center}
\begin{tabular}{|c|c|c|c|}
\hline 
Values of $n$ such that $a(n)=0$ & $\frac{(n+1)^{3}-1}{n}$ & $\delta(n)$ & $\delta(\frac{(n+1)^{3}-1}{n})$\tabularnewline
\hline
\hline 
2 & 23 & 1 & 1\tabularnewline
\hline 
23 & 601 & 1 & 1\tabularnewline
\hline 
119333 & 142432291 & 1 & 1\tabularnewline
\hline
\end{tabular}
\par\end{center}

We can also start from another initial value. For instance let $p=2$
and choose $a(1)=2$ this gives

\begin{center}
\begin{tabular}{|c|c|c|c|}
\hline 
Values of $n$ such that $a(n)=0$ & $n+2$ & $\delta(n)$ & $\delta(n+2)$\tabularnewline
\hline
\hline 
3 & 5 & 1 & 1\tabularnewline
\hline 
17 & 19 & 1 & 1\tabularnewline
\hline 
281 & 283 & 1 & 1\tabularnewline
\hline 
79559 & 79561 & 1 & 1\tabularnewline
\hline 
6329815697 & 6329815699 & 1 & 1\tabularnewline
\hline
\end{tabular}
\par\end{center}

See the APPENDIX 3 for more experiments supporting the conjecture
for $p=2$ and various starting values.

\subsection{Conjecture 3bis}

In the same vein we find an algorithm seeming generating primes pairs
of type $(p,p+2m)$. Let $a(1)=4m^{2}$ and define:
\begin{itemize}
\item $a(n)=\left|a(n-1)-\gcd(a(n-1),n(n+2m))\right|$ 
\end{itemize}
Then we claim that for $n$ large enough:
\begin{itemize}
\item $a(n)=0\Rightarrow$$n+1$ and $n+2m+1$ are primes and this happens
infinitely many times and there are infinitely many pairs $(n+1,n+2m+1)$
of consecutive primes.
\end{itemize}
Thus the Polignac conjecture would be true. Here a table for $2m=4$

\begin{center}
\begin{tabular}{|c|c|c|c|}
\hline 
Values of $n$ when $a(n)=0$ & $n+5$ & $\delta(n+1)$ & $\delta(n+5)$\tabularnewline
\hline
\hline 
12 & 17 & 1 & 1\tabularnewline
\hline 
192 & 197 & 1 & 1\tabularnewline
\hline 
38196 & 38201 & 1 & 1\tabularnewline
\hline 
1459118862 & 1459118867 & 1 & 1\tabularnewline
\hline
\end{tabular}
\par\end{center}

See APPENDIX 4 for experiments with other small values of $m$. The
conjecture 9 will give an easier way for computation of such pairs
of primes since we avoid the square in the gcd.

\paragraph*{Remark}

It is also possible to use Rowland-Shevelev recursion for generating
increasing sequences of twin primes or things like that. For instance
let:
\begin{itemize}
\item $a_{1}=2$ and $a_{n}=a_{n-1}+\gcd(a_{n-1},n(n-2))$ 
\end{itemize}
Then the sequence of records for the differences $a_{n}-a_{n-1}$
yields an increasing sequence of lower of twin primes. See the end
of the APPENDIX 3.

\subsection{Conjecture 3ter}

We generate primes triplet of type $(p,p+2,p+6)$. Let $a(1)=4$ and
define:
\begin{itemize}
\item $a(n)=\left|a(n-1)-\gcd(a(n-1),n(n+2)(n+6))\right|$ 
\end{itemize}
Then we claim 
\begin{itemize}
\item $a(n)=0\Rightarrow$$(n+1,n+3,n+7)$ is a prime triplet.
\end{itemize}
Here a table 

\begin{center}
\begin{tabular}{|c|c|c|c|}
\hline 
Values of $n$ such that $a(n)=0$ & $\delta(n+1)$ & $\delta(n+3)$ & $\delta(n+7)$\tabularnewline
\hline
\hline 
4 & 1 & 1 & 1\tabularnewline
\hline 
40 & 1 & 1 & 1\tabularnewline
\hline 
82006 & 1 & 1 & 1\tabularnewline
\hline
\end{tabular}
\par\end{center}

Although it is hard to perform convincing experiments, it is coherent
with our previous conjectures and some very general rule should exist.
This conjecture can be extended to any sort of prime triplet and to
$m$-uplet. The conjecture 9 is better fitted for practical computation
of $m$-uplet.

\section{Finding big primes}

Here we merge the conjectures 2 and 3 to obtain another conjectural
way to unearth increasing sequences of primes. The rate of growth
is multiple exponential and this could produce very big primes compared
to the number of required iterations. We then discuss probabilities
issues and propose a method to got very big primes similarly as what
is described in 1.3.

\subsection{Conjecture 4}

Suppose $p\geq2$ is prime and let $a(1)=p$ and:
\begin{itemize}
\item $a(n)=\left|a(n-1)-\gcd(a(n-1),pn^{2}-1)\right|$
\end{itemize}
Then we claim:
\begin{itemize}
\item there are infinitely many values of $n$ such that $a(n)=0$.
\item for $n$ large enough (usually $n>2$ is working for small $p$) we
have $a(n)=0\Rightarrow$$p(n+1)^{2}-1$ is prime.
\end{itemize}
See APPENDIX 5 for the begining of tables for some values of $p$.
Could we by chance get so many prime numbers? 

$ $

\subsection{Why stopping here?}

In fact it seems one can go further. Suppose $a(1)=2$ and:
\begin{itemize}
\item $a(n)=\left|a(n-1)-\gcd(a(n-1),2n^{3}-1)\right|$
\end{itemize}
Then we claim:
\begin{itemize}
\item there are infinitely many values of $n$ such that $a(n)=0$.
\item $a(n)=0\Rightarrow$$2(n+1)^{3}-1$ is prime for $n$ large enough.
\end{itemize}
\begin{flushleft}
Here the begining of the table
\par\end{flushleft}

\begin{center}
\begin{tabular}{|c|c|c|}
\hline 
Values of $n$ such that $a(n)=0$ & $2(n+1)^{3}-1$ & $\delta(2(n+1)^{3}-1)$\tabularnewline
\hline
\hline 
3 & 127 & 1\tabularnewline
\hline 
125 & 4000751 & 1\tabularnewline
\hline 
4000877 & 128084306502569672303  & 1\tabularnewline
\hline
\end{tabular}
\par\end{center}

However it is not easy to generalise this one. We can also find good
initial values and/or a good factor before $n^{3}$. For instance
let $w(1)=3$ and:
\begin{itemize}
\item $a(n)=\left|a(n-1)-\gcd(a(n-1),10n^{3}-1)\right|$
\end{itemize}
Then we get

\begin{center}
\begin{tabular}{|c|c|c|}
\hline 
Values of $n$ such that $a(n)=0$ & $10(n+1)^{3}-1$ & $\delta(10(n+1)^{3}-1)$\tabularnewline
\hline
\hline 
4 & 1249 & 1\tabularnewline
\hline 
1240 & 19112405209 & 1\tabularnewline
\hline
\end{tabular}
\par\end{center}

And we can provide this other impressive example with exponent 7.
Suppose $a(1)=3$ and:
\begin{itemize}
\item $a(n)=\left|a(n-1)-\gcd(a(n-1),2n^{7}-1)\right|$
\end{itemize}
\begin{center}
\begin{tabular}{|c|c|c|}
\hline 
Values of $n$ such that $a(n)=0$ & $2(n+1)^{7}-1$ & $\delta(2(n+1)^{7}-1)$\tabularnewline
\hline
\hline 
2 & 4373 & 1\tabularnewline
\hline 
4352 & 59231218330987879606185473 & 1\tabularnewline
\hline
\end{tabular}
\par\end{center}

\subsection{How many chances we have to catch a prime?}

Continuing this way it seems possible to find very big primes using
the $(k,b,c)$ recursion:
\begin{itemize}
\item $a(1)=k$ and $a(n)=\left|a(n-1)-\gcd(a(n-1),bn^{c}-1)\right|$ 
\end{itemize}
with suitable choices of $(k,b,c).$ A natural question is then: what
is the chance to get a prime when the algorithm reachs the first zero
starting with any value $N$? To evaluate this chance let $r_{k}=\min\left\{ i\geq1\mid a(i)=0\right\} $
and:
\begin{itemize}
\item $\Upsilon(N)=\frac{1}{N}\#\left\{ k\mid1\leq k\leq N\,\&\,\delta(b(r_{k}+1)^{c}-1)=1\right\} $ 
\end{itemize}
wich represents the chance to get a prime reached by the algorithm
starting with $N$. As $N\rightarrow\infty$ it appears this chance
is not zero. For instance if $(b,c)=(2,2)$ we have $\Upsilon(N)\simeq0.8$
as $N\rightarrow\infty$. Which gives an efficient method to got big
primes since the first zero is reached after a number of iterations
of order $N$ and thus we have more than $80\%$ of chance to got
a prime of size $2N^{2}$. Here is a graph supporting this claim.
We plot $\Upsilon(N)$ for $N=1,2,3,\ldots,20000$ and for $(b,c)=(3,2)$
(pink) $(b,c)=(5,2)$ (blue)

\begin{center}
(fig.2)
\par\end{center}

\begin{center}
\includegraphics[width=0.5\paperwidth]{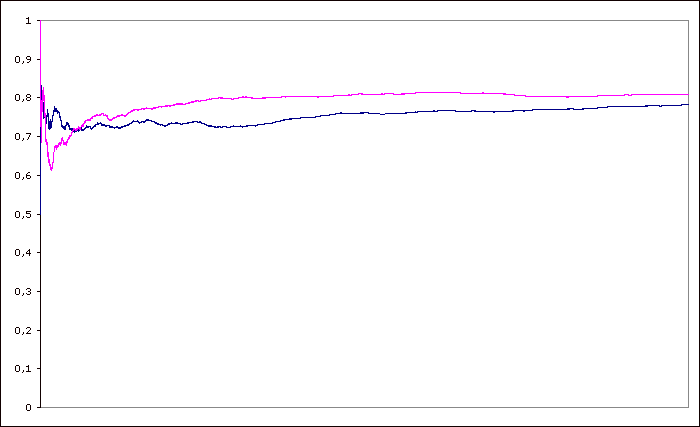}
\par\end{center}

And we believe $\lim_{N\rightarrow\infty}\Upsilon(N)=L(c)$ exists
and depends only on $c$ and it is clear that $c<c'\Rightarrow L(c)>L(c')$.

\subsection{A rule of construction}

As seen before in 1.3. there is here also a general relation between
consecutive records. Namely we still consider for any initial value:
\begin{itemize}
\item $a(n)=\left|a(n-1)-\gcd(a(n-1),bn^{c}-1)\right|$ 
\end{itemize}
and we define $\left(w(k)\right)_{k\geq1}$ as the sequence of values
taken by $a_{n}-a_{n-1}$ when $\left|a_{n}-a_{n-1}\right|>1$. Then
the records values of $\left(w(k)\right)_{k\geq1}$ are given by $k$
such that $w(k)>0$. Let us now write this increasing sequence of
$k$ using a sequence $\left(\alpha_{j}\right)_{j\geq1}$. Then we
have the following simple relationship between 2 consecutive records
values $w(\alpha_{j})$ and $w(\alpha_{j+1})$ (details ommited):
\begin{itemize}
\item $w(\alpha_{j+1})=b\left(w(\alpha_{j})+1+\left(\frac{w(\alpha_{j})+1}{b}\right)^{1/c}+\sum_{i=\alpha_{j}+1}^{\alpha_{j+1}-1}\left(a_{i+1}-a_{i}+1\right)\right)^{c}-1$ 
\end{itemize}
Thus as in 1.3.1. experiments show that it isn't necessary to compute
all terms in the sum to got the next record value. Indeed there is
again a value $n_{0}$ conjectured to be closed of $\alpha_{j}$ such
that $\alpha_{j}<n_{0}\leq i<\alpha_{j+1}\Rightarrow a_{i+1}-a_{i}=-1$.
Therefore this gives sometime an efficient method for building a bigger
prime from a prime record value or a non prime record value. The quadratic
case seems well working. For instance let us consider this quadratic
case:
\begin{itemize}
\item $a(1)=1$ and $a(n)=\left|a(n-1)-\gcd(a(n-1),32n^{2}-1)\right|$ 
\end{itemize}
Then we get the sequence of values $a_{n}-a_{n-1}$ for those $n$
$\leq10^{6}$ such that $\left|a_{n}-a_{n-1}\right|>1$:
\begin{lyxlist}{00.00.0000}
\item [{$127,-7,-17,-7,-7,294911,-1289,2760686028799,-113,-103,-7,-7,-113.$}]~
\end{lyxlist}
The 3 first records $\left(127,294911,2760686028799\right)$are prime
values and supposing there is no more value until the next record
(say $X$) the formula above yields:
\begin{lyxlist}{00.00.0000}
\item [{$X=32\left(2760686028800+\left(\frac{2760686028800}{32}\right)^{1/2}-112-102-6-6-112\right)^{2}-1$}]~
\end{lyxlist}
Giving $X=243884447023448880167715967$ which is still a prime value.
We provide also an exemple for the cubic case:
\begin{itemize}
\item $a(1)=2$ and $a(n)=\left|a(n-1)-\gcd(a(n-1),2n^{3}-1)\right|$ 
\end{itemize}
Then we get the sequence of values $a_{n}-a_{n-1}$ for those $n$
$\leq10^{5}$ such that $\left|a_{n}-a_{n-1}\right|>1$:
\begin{lyxlist}{00.00.0000}
\item [{$-3,53,-5,-3,265301,-109,-31,-17,-3,-5,-3.$}]~
\end{lyxlist}
So we suspect there is no more value until the next record value (say
$X$) and the formula above yields:
\begin{lyxlist}{00.00.0000}
\item [{$X=2\left(265301+1+\left(\frac{265301+1}{2}\right)^{1/3}-108-16-2-4-2\right)^{3}-1$}]~
\end{lyxlist}
giving $X=37299785868725741$ which indeed is the next record value
and is a prime value. 

For higher exponent it seems less easy to find many working exemples
but we think it would be worth to explore this method further in order
to check its possible efficiency. The main question would be: are
there any conditions forcing records values to be prime values and
making $n_{0}$ very close from the working record value?

\subsection{Generating big twin primes}

We can also do the same kind of task (cf. 3.3.) for twin primes borrowing
from Shevelev the idea for a quadratic case. Let:
\begin{itemize}
\item $a(0)=k$ and $a(n)=a(n-1)-\gcd(a(n-1),2n^{2}+(-1)^{n})$ . 
\end{itemize}
Define:
\begin{itemize}
\item $\Upsilon_{twin}(N)=\frac{1}{N}\#\left\{ k\mid1\leq k\leq N\,\&\,\delta(2(r_{k}+1)^{2}-1)\delta(2(r_{k}+1)^{2}+1)=1\right\} $ 
\end{itemize}
wich represents the chance to get a pair of twin primes reached by
the algorithm starting with $N$. We plot below $\Upsilon_{twin}(N)$
for $N=1,2,3,\ldots,20000$

\begin{center}
(fig.3)
\par\end{center}

\begin{center}
\includegraphics[width=0.5\paperwidth]{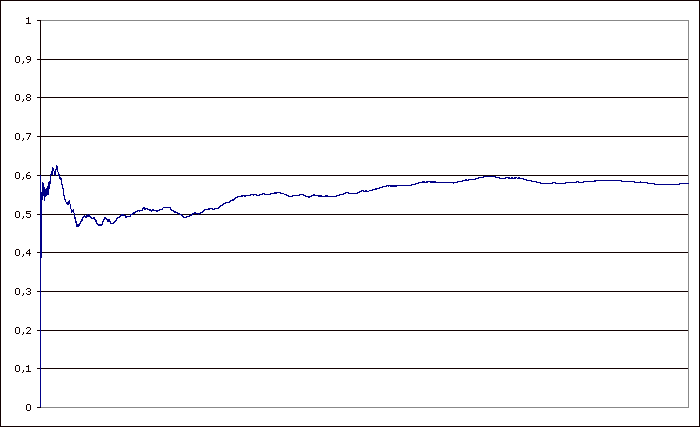}
\par\end{center}

And the graph in the range $2\leq N\leq30000$ of $\frac{r_{N}}{N}\delta(2(r_{N}+1)^{2}-1)\delta(2(r_{N}+1)^{2}+1)$ 

\begin{center}
(fig.4)
\par\end{center}

\begin{center}
\includegraphics[width=0.5\paperwidth]{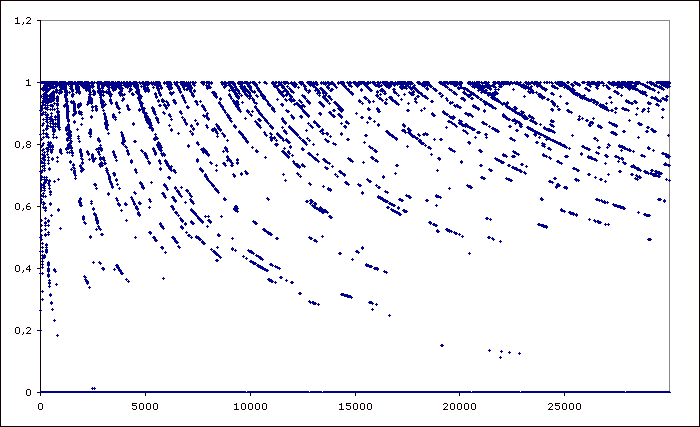}
\par\end{center}

The fact the $x$ axis is blue means there are many zero values.

Roughly speaking it shows we have $50\%$ chance to come across a
twin prime of size $N^{2}$ after $N$ iterations starting from a
large random value of $N$. Certainly one should search for simple
conditions on $N$ in order to increase this chance. For instance
let us start from $2N^{2}$ when $2N^{2}-1$ and $2N^{2}+1$ are already
primes. This sequence of $N$ begins:
\begin{itemize}
\item $3,6,21,24,36,42,45,87,102,132,153,186,204,228,237,273,297,300,321,...$
\end{itemize}
Let now $v(n)$ denotes the $n$-th term of this sequence and let:
\begin{itemize}
\item $\Upsilon_{v}(N)=\frac{1}{N}\#\left\{ k\mid1\leq k\leq N\,\&\,\delta(2(r_{v(k)}+1)^{2}-1)\delta(2(r_{v(k)}+1)^{2}+1)=1\right\} $
\end{itemize}
We plot below $\Upsilon_{v}(N)$ for $N=1,2,3,\ldots,140$

\begin{center}
(fig.5)
\par\end{center}

\begin{center}
\includegraphics[width=0.5\paperwidth]{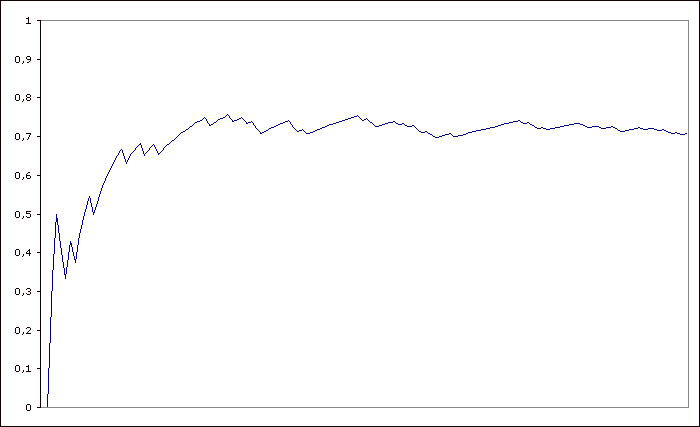}
\par\end{center}

So it seems we have slightly more chances to got a twin prime pair
starting with $2v(n)^{2}$ ($\sim70\%$ of chance to got a twin prime
of order $n^{4}$) than starting with $n$ ($\sim50\%$ of chance
producing primes of order $n^{2}$). Although some computation is
needed to got the sequence $v(n)$ it could be an efficient method
for generating big primes. More importantly, this observation is a
striking one regarding the conjectures like 2 or 3. Indeed, it appears
that when we catch a prime with a given property and run the gcd-algortihm
another time (starting around this value) we have more chance to get
a new prime with this property. Heuristically this should explain
why there are apparently infinite {}``chains'' of primes generated
by the algorithm in conjectures 2, 3 or similar ones when we have
a good starting value. Repeating the process seems to force the algorithm
to reach primes every time once we are on a right track.

\section{Conjecture 5 : on the Goldbach conjecture}

We propose a first constructive way to prove this famous conjecture.
This is somewhat unsatisfactory since there are very few non working
starting values in the range where we performed the computation. Consider
$N\geq2$ and let $a(1)=N-2$ $ $ and define for $n\geq2$:
\begin{itemize}
\item $a(n)=a(n-1)-\gcd(a(n-1),(n-1)(2N-n+1))$ 
\end{itemize}
Then we claim there is always a unique $g_{N}\in\left\{ 2,3,...,N-2\right\} $
such that :
\begin{itemize}
\item $a(g_{N})=0$ 
\end{itemize}
and we have:
\begin{itemize}
\item $g_{N}$ and $2N-g_{N}$ are simultanuously primes except for very
few $N$ ( a set conjectured to be of measure zero). 
\end{itemize}
Hereafter we plot $\frac{g_{N}}{N}\delta(g_{N})\delta(2N-g_{N})$
for $2\leq N\leq30000$. 

\begin{center}
(fig.6)
\par\end{center}

\begin{center}
\includegraphics[width=0.5\paperwidth]{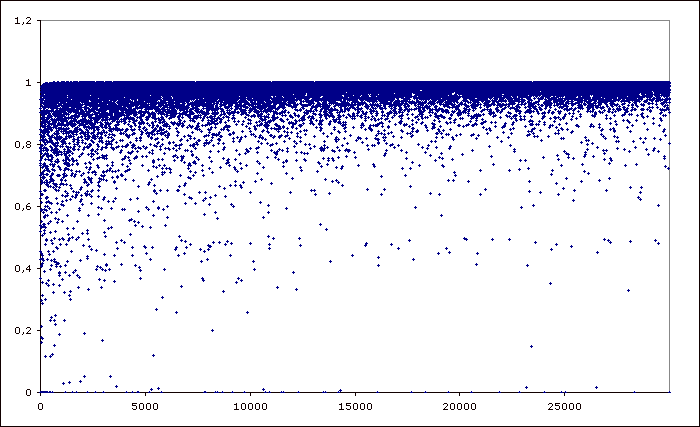}
\par\end{center}

We can see on the $x$ axis some zero values which become very sparse
when $N$ increases (perhaps there is no more zero value for $N$
sufficiently large). In the conjecture 10 we provide a variation of
this conjecture where clearly there is no exceptional set of non working
values for $N$ large enough.

\section{The Schinzel's hypothesis H }

Finally this kind of method should have considerable application.
One of them, regarding what we discuss before, could be a new version
of the Schinzel's hypothesis H \cite{key-5}. Indeed, suppose $P(x)$
is a polynomial with integer coefficients and define the sequence
$S$ as follows. 
\begin{itemize}
\item $S(1)\in\mathbb{N}$
\item $S(n)=\left|S(n-1)-\gcd(S(n-1),P(n))\right|$ 
\end{itemize}
Under some assumptions the values of $n$ such that $S(n)=0$ imply
$P(n+1)$ has nice arithmetical properties. This is what we have seen
previously. But let us consider three additional examples.

\subsection*{$P(x)=x^{2}+1$}

We claim that for a suitable starting value $S(1)$ such as $S(1)=2$
we have:
\begin{itemize}
\item $S(n)=0$ for infinitely many values of $n$.
\end{itemize}
$S(n)=0\Rightarrow P(n+1)$ is a prime number.

Here a table with the first values

\begin{center}
\begin{tabular}{|c|c|c|}
\hline 
Values of $n$ such that $S(n)=0$ & $P(n+1)$ & $\delta(P(n+1))$\tabularnewline
\hline
\hline 
3 & 17 & 1\tabularnewline
\hline 
13 & 197 & 1\tabularnewline
\hline 
203 & 41617 & 1\tabularnewline
\hline 
41813 & 1748410597  & 1\tabularnewline
\hline
\end{tabular}
\par\end{center}

\subsection*{$P(x)=x^{3}+1=(x+1)(x^{2}-x+1)$}

We claim that for a suitable starting value $S(1)$ such as $S(1)=2$
we have:
\begin{itemize}
\item $S(n)=0$ for infinitely many values of $n$.
\item $S(n)=0\Rightarrow(n+2)$ and $(n^{2}+n+1)$ are prime numbers.
\end{itemize}
Here a table with the first values

\begin{center}
\begin{tabular}{|c|c|c|c|}
\hline 
Values of $n$ such that $S(n)=0$ & $n+2$ & $n^{2}+n+1$ & $\delta(n+2)$$\delta(n^{2}+n+1)$\tabularnewline
\hline
\hline 
3 & 5 & 13 & 1\tabularnewline
\hline 
69 & 71 & 4831 & 1\tabularnewline
\hline 
299391 & 299393 & 89635270273 & 1\tabularnewline
\hline
\end{tabular}
\par\end{center}

\subsection*{$P(x)=(2x-3)(x^{2}-x+1)$}

We claim that for a suitable starting value $S(1)$ we have:
\begin{itemize}
\item $S(n)=0$ for infinitely many values of $n$.
\item $S(n)=0\Rightarrow(2n-1)$ and $(n^{2}+n+1)$ are prime numbers.
\end{itemize}
Here a table with the first values for $S(1)=1,2,5$ 

\begin{center}
\begin{tabular}{|c|c|c|c|}
\hline 
$S(1)$ & $n$ such that $S(n)=0$ & $2n-1$ & $n^{2}+n+1$\tabularnewline
\hline
\hline 
1 & 2 & 3 & 7\tabularnewline
\hline 
 & 24 & 47 & 601\tabularnewline
\hline 
 & 24186 & 48371 & 584986783\tabularnewline
\hline 
2 & 3 & 5 & 13\tabularnewline
\hline 
 & 69 & 137 & 4831\tabularnewline
\hline 
 & 658657 & 1317713 & 434093205307\tabularnewline
\hline 
5 & 6 & 11 & 43\tabularnewline
\hline 
 & 414 & 827 & 171811\tabularnewline
\hline 
 & 141629682 & 283259363 & 20058966965050807\tabularnewline
\hline
\end{tabular}
\par\end{center}

And all the values for $2n-1$ and $n^{2}+n+1$ in the table are primes. 

Above exemples allow us to unify all our previous observations (except
for the conjecture $5$ which is a special case needing more thought
to be generalized). This is the next conjecture.

\subsection{Conjecture 6}

Clearly something very general is working and one can imagine to state
a deep conjecture. Suppose $P(x)=\prod_{j=1}^{m}Q_{j}(x)$ where $P$
is polynomial with integer coefficients (not all zeros). Suppose each
$Q_{j}$ is irreducible and $Q_{1}(k),Q_{2}(k),\ldots,Q_{m}(k)$ can
simultanously be primes for large $k$. Then we claim there are infinitely
many values of $S(1)$ such that:
\begin{itemize}
\item $S(n)=0$ for infinitely many values of $n$.
\item $S(n)=0\Rightarrow Q_{j}(n+1)$ is simultaneously prime for $1\leq j\leq m$.
\end{itemize}

\subsubsection{Remark}

This could have direct application such as generating not only big
primes (see conjecture 4) but also big primes with given property
like twin primes. For instance let:
\begin{itemize}
\item $P(x+1)=(x{}^{2}+1)(x^{2}+3)$
\item $S(1)=4$
\end{itemize}
Then we compute big twin primes compared to the number of iterations:

\begin{center}
\begin{tabular}{|c|c|c|c|}
\hline 
Values of $n$ such that $S(n)=0$ & $n^{2}+1$ & $n^{2}+3$ & $\delta(n^{2}+1)$$\delta(n^{2}+3)$\tabularnewline
\hline
\hline 
2 & 5 & 7 & 1\tabularnewline
\hline 
14 & 197 & 199 & 1\tabularnewline
\hline 
32374 & 1048075877 & 1048075879 & 1\tabularnewline
\hline
\end{tabular}
\par\end{center}

\subsection{Conjecture 7}

Instead of considering an infinite sequence we consider a starting
value and see what happens. This gives a modified version of the conjecture
6 with more details suggesting some analytic study.\textcompwordmark{}

Suppose that $P(x)$ is a polynomial with integer coefficients satisfying
$P(n)\geq0$ and $P(+\infty)=+\infty$. Let us define the sequence
$a$ as follows. 
\begin{itemize}
\item $a(1)=N\in\mathbb{N}$
\item $a(n)=a(n-1)-\gcd(a(n-1),P(n))$ 
\end{itemize}
Then we conjecture without any other assumption on $P$:
\begin{itemize}
\item $\exists\, f(N)\in\left\{ 1,2,\ldots,N\right\} $ such that $a(f(N))=0$.
\item $\liminf_{N\rightarrow\infty}\frac{f(N)}{N}\geq0$ and $\limsup{}_{N\rightarrow\infty}\frac{f(N)}{N}=1$.
\end{itemize}
Now suppose as above $P(x)=\prod_{j=1}^{m}Q_{j}(x)$ where $P$ is
polynomial with integer coefficients (not all zeros). Suppose each
$Q_{j}$ is irreducible and $Q_{1}(k),Q_{2}(k),\ldots,Q_{m}(k)$ can
be simultanously prime for large $k$. Then there are infinitely many
$N$ such that we get:
\begin{itemize}
\item $Q_{j}(f(N)+1)$ is prime for $j\in\left\{ 1,2,\ldots,m\right\} $
\end{itemize}
And more precisely there is a positive proportion of $N$ such that
$Q_{j}(f(N)+1)$ is prime for $j\in\left\{ 1,2,\ldots,m\right\} .$
Let us define this proportion as follows:
\begin{itemize}
\item $L(P)=\lim_{n\rightarrow\infty}\frac{1}{n}\#\left\{ k\mid1\leq k\leq n\,\&\,\prod_{j=1}^{m}\delta(Q_{j}(f(k)+1))=1\right\} $.
\end{itemize}
Then we claim:
\begin{itemize}
\item $Q_{j}$ is of degree $1$ for any $j\in\left\{ 1,2,\ldots,m\right\} \Rightarrow L(P)=1$. 
\item One $Q_{j}$ is irreducible of degree 2 $\Rightarrow L(P)<1$ (see
3.3. with a sample giving $L(2x^{2}-1)\sim0.8$). 
\end{itemize}
Here the graph of $\frac{f(N)}{N}\delta(f(N)+1)$ when $P(x)=x$ for
$1\leq N\leq20000$ 

\begin{center}
(fig.7)
\par\end{center}

\begin{center}
\includegraphics[width=0.5\paperwidth]{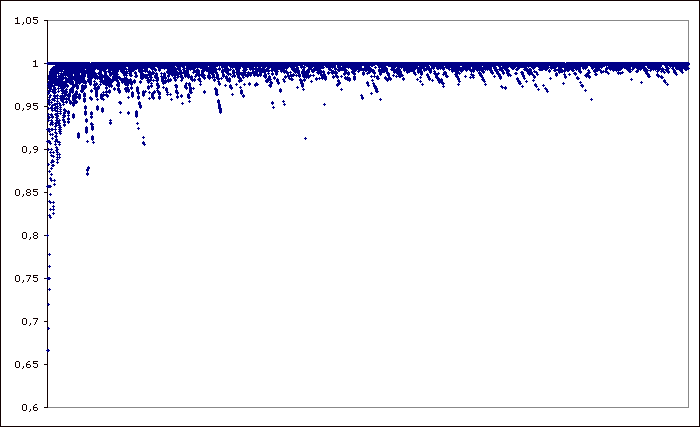}
\par\end{center}

Here the graph of $\frac{f(N)}{N}\delta(2f(N)+3)\delta(f(N)+3)\delta(f(N)+7)$
when $P(x)=(2x+1)(x+2)(x+6)$ for $N=100k$ and $1\leq k\leq20000$ 

\begin{center}
(fig.8)
\par\end{center}

\begin{center}
\includegraphics[width=0.5\paperwidth]{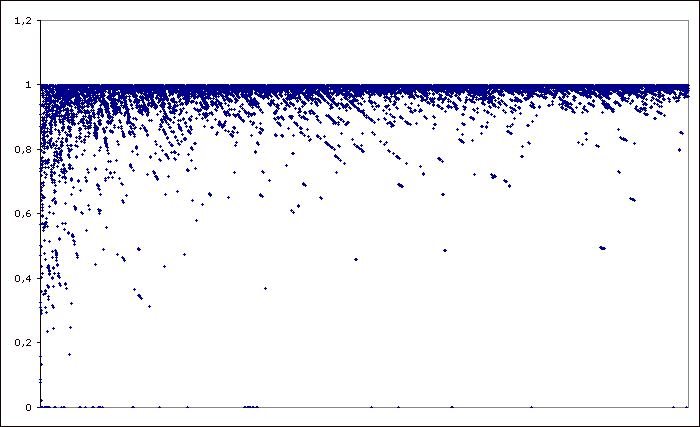}
\par\end{center}

One can see the values of $N$ such that $\delta(2f(N)+3)\delta(f(N)+3)\delta(f(\alpha N)+7)=0$
become sparse when $N$ increases. We believe that for $N$ large
enough this set is of measure zero since $Q_{j}$ are of degree $1$.

\section{Conjecture 8}

We propose a way to prove again there are infinitely many twin primes
but in a slightly different way than before. The recursion is the
backwards version of the Shevelev recursion. First we give a conjectural
way to prove there are infinitely many prime numbers (this is very
similar to what we said in conjecture 5 and should be easy to prove).
Although it is not the simplest way to prove there are infinitely
many primes, it seems important to start with the usual primes in
order to see better how it is also working for twin primes. The two
cases are indeed very similar.

\subsection{There are infinitely many primes}

Let:
\begin{itemize}
\item $a(1)=N-2\geq0$ 
\item $n\geq2\Rightarrow$$a(n)=a(n-1)-\gcd(a(n-1),n-1)$ 
\end{itemize}
Then we claim we have:
\begin{enumerate}
\item $\forall N\geq4$, $\exists f(N)\in\left\{ 2,\ldots,N-1\right\} $
such that $a(f(N))=0$.
\item $f(N)\sim N$ $\left(N\rightarrow\infty\right)$ and more precisely
we claim $f(N)=N+o(N^{1/2}\log N)$
\item $f(N)$ is prime.
\item $f(N)=N-1\Leftrightarrow N-1$ is a prime number.
\item $f(N)=N-2\Leftrightarrow N-2$ is an odd prime number.
\item $f(N)=N-3$ or $N-4$ $\Leftrightarrow N-3$ (resp $N-4$) is a prime
number of form $6k+1$.
\item $f(N)=N-5$ or $N-6$ $\Leftrightarrow N-5$ (resp $N-6$) is a prime
number of form $30k+1$.
\end{enumerate}
Since $f(N)\rightarrow\infty$ as $N\rightarrow\infty$ (from 2.)
there are infinitely many primes (from 3.). 

Here the graph of $\frac{f(N)}{N-1}\delta(f(N))$ for $2\leq N\leq20000$.

\begin{center}
(fig.9)
\par\end{center}

\begin{center}
\includegraphics[width=0.5\paperwidth]{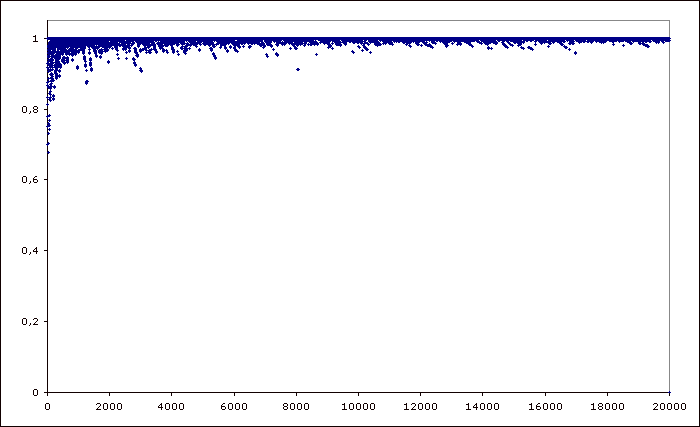}
\par\end{center}

To see the more precise behaviour claimed in 2. we plot $\frac{N-f(N)}{n^{1/2}}$
for $3\leq N\leq20000$.

\begin{center}
(fig.10)
\par\end{center}

\begin{center}
\includegraphics[width=0.5\paperwidth]{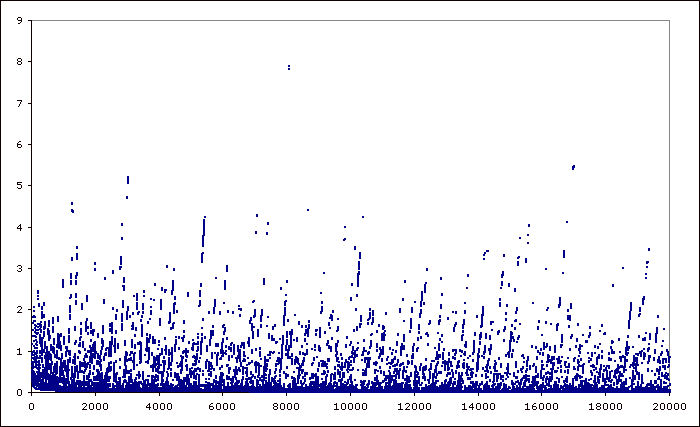}
\par\end{center}

This is quite erratic and divided by $\log n$ it should converge
very slowly to zero but more experiments are needed to confirm the
tendancy.

\subsubsection*{Remark}

It is worth to compare this sequence $f(N)$ to $p(\pi(N))$ where
$p(n)$ denotes the $n$-th prime and $\pi(x)$ is the prime couting
function. So that $p(\pi(N))$ is the largest prime $\leq N$. Here
we still plot $\frac{N-p(\pi(N))}{n^{1/2}}$ for $3\leq N\leq20000$.

\begin{center}
(fig.11)
\par\end{center}

\begin{center}
\includegraphics[width=0.5\paperwidth]{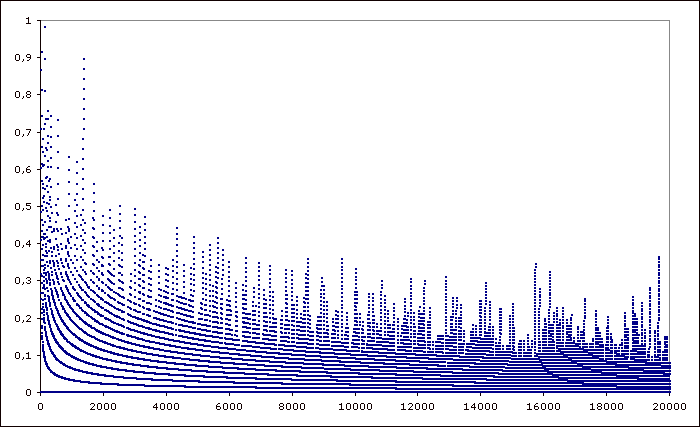}
\par\end{center}

Things are less erratic and here the graph goes more certainly to
zero without dividing by $\log n$. We will discuss about an important
consequence of this conjectured behaviour in section 9.

\subsection{There are infinitely many twin primes}

Let:
\begin{itemize}
\item $a(1)=N-2\geq0$ 
\item $a(n)=a(n-1)-\gcd(a(n-1),n+(-1)^{n})$ 
\end{itemize}
Then we claim we have:
\begin{enumerate}
\item $\forall N\geq2$, $\exists h(N)\in\left\{ 1,2,\ldots,N\right\} $
such that $a(h(N))=0$.
\item $h(N)\sim N$ $\left(N\rightarrow\infty\right)$.
\item $\forall N\geq98$, $h(N)$ and $h(N)+2$ are necessarily simultanously
primes.
\item For $N\geq4$ we have $h(N)=N-1\Leftrightarrow(N-1,N+1)$ is a pair
of twin primes.
\item For $N\geq13$ we have $h(N)=N-2\Leftrightarrow N$ is the greater
of a pair of twin primes.
\item For $N\geq14$ we have $h(N)=N-3\Leftrightarrow N-3$ is the lesser
of a pair of twin primes.
\end{enumerate}
Since $h(N)\rightarrow\infty$ as $N\rightarrow\infty$ (from 2) there
are infinitely many twin primes (from 3).

Here the graph of $\frac{h(N)}{N}\delta(h(N))\delta(h(N)+2)$ for
$1\leq N\leq30000$

\begin{center}
(fig.12)
\par\end{center}

\begin{center}
\includegraphics[width=0.5\paperwidth]{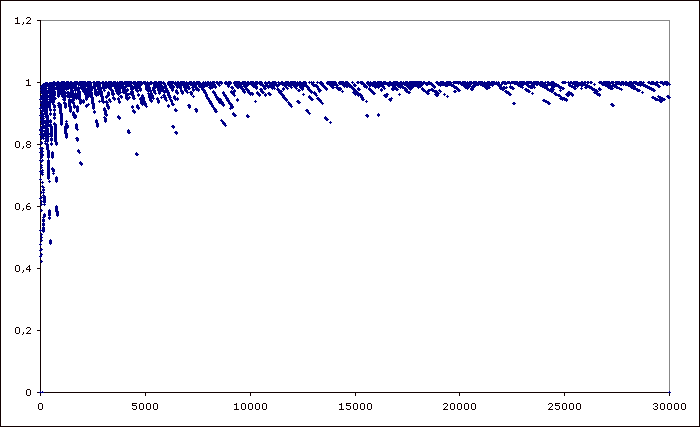}
\par\end{center}

We see there is no more zero value for $N\geq98$.

\section{Conjecture 9}

Let $m\in\mathbb{N}$ and define the sequence $a$ as follows:
\begin{itemize}
\item $a(1)\in\mathbb{N}$ 
\item $a(n)=\left|a(n-1)-\gcd(a(n-1),mn+b(n))\right|$ 
\end{itemize}
Where $(b_{n})_{n\geq1}$ is a periodic sequence of period length
$\beta$ such that:
\begin{itemize}
\item $mn+b_{1},mn+b_{2},\ldots,mn+b_{\beta}$ can be simultanously prime
for $n$ large.
\end{itemize}
Then we claim we have:
\begin{itemize}
\item $a(n)=0$ for infinitely many values of $n$.
\item for $n$ large enough we have $a(n)=0\Rightarrow$ $m(n+1)+b_{1},m(n+1)+b_{2},\ldots,m(n+1)+b_{\beta}$
are simultanously prime.
\item the rate of growth of the sequence of $n_{i}$ such that $a(n_{i})=0$
is like $(m+1)^{i}$ . 
\end{itemize}
See APPENDIX 6 for experiments where we provide examples for $\beta$-periodic
sequences and $\beta\in\left\{ 2,3,4,5,6\right\} $.

\subsection*{Generalisation}

Although we check few cases of this generalisation it is worth to
mention this result which is an interesting variation of the conjectures
6 and 9. Suppose $Q_{j}$ is irreducible and $Q_{1}(k),Q_{2}(k),\ldots,Q_{\beta}(k)$
can be simultanously prime for large $k$. Suppose $(b_{n})_{n\geq1}$
is a periodic sequence of period length $\beta$ such that $\left\{ b(i)\right\} _{1\leq i\leq\beta}$
is a permutation of $\left\{ i\right\} _{1\leq i\leq\beta}$ and define
the sequence $a$ as follows:
\begin{itemize}
\item $a(1)\in\mathbb{N}$ 
\item $a(n)=\left|a(n-1)-\gcd(a(n-1),Q_{b(n)}(n))\right|$ 
\end{itemize}
Then we claim we have:
\begin{itemize}
\item $a(n)=0$ for infinitely many values of $n$.
\item for $n$ large enough we have $a(n)=0\Rightarrow$ $Q_{1}(n+1),Q_{2}(n+1),\ldots,Q_{\beta}(n+1)$
are simultanously prime.
\end{itemize}

\subsection*{Extension of the conjecture 6}

We can extend what we done in conjecture 6 (and also in conjecture
7) with twin primes to prime triplet or any sort of $m$-uplet. For
instance let us see how this is also working with prime triplet of
type $(p,p+2,p+6)$. Let $\left(b(n)\right)_{n\geq1}$ be the 3-periodic
sequence $\left\{ 2,6,0\right\} $ and define the sequence:
\begin{itemize}
\item $a(n)=N-6\geq0$ 
\item $a(n)=a(n-1)-\gcd(a(n-1),n+b(n))$ 
\end{itemize}
Then we claim:
\begin{enumerate}
\item $\forall N\geq6$ $\exists f(N)\in\left\{ 1,2,\ldots,N\right\} $
such that $a(f(N))=0$.
\item $f(N)\sim N$ $\left(N\rightarrow\infty\right)$.
\item $\forall N\geq2735$, $\left(f(N)+1,f(N)+3,f(N)+7\right)$ is a prime
triplet.
\item For $N\geq5$ we have $h(N)=N-1\Leftrightarrow(N,N+2,N+6)$ is a prime
triplet of type $(p,p+2,p+6)$ (sequence $A022004$ in \cite{key-2}
)
\end{enumerate}
Here the graph of $\frac{f(N)}{N}\delta(f(N)+1)\delta(f(N)+3)\delta(f(N)+7)$
for $6\leq N\leq30000$

\begin{center}
(fig.13)
\par\end{center}

\begin{center}
\includegraphics[width=0.5\paperwidth]{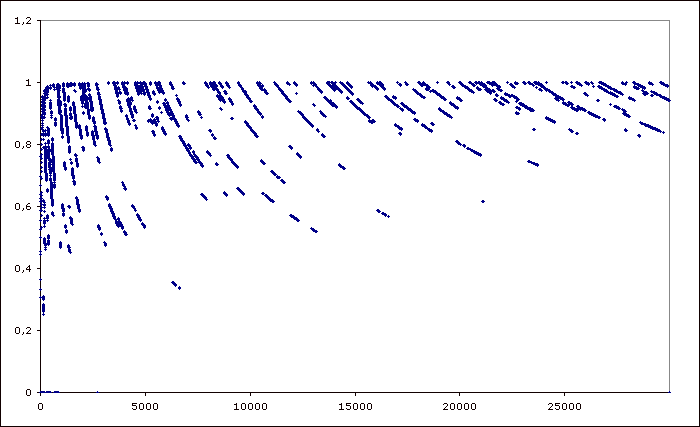}
\par\end{center}

We see that for $N\geq2735$ there is no more zero value.

\section{Conjecture 10: on the Goldbach conjecture }

Finally we found a way to adapt the method to the Goldbach conjecture.
This was not easy since this recursion is very sensitive (to the initial
value and to what we put in the gcd). This formulation of the conjecture
is nicer than the conjecture 5 and is very similar to the conjecture
6 and our formulation of the weak twin prime conjecture. Indeed the
method works in both cases for $N$ large enough. Namely define the
sequence $a$ for $N\geq2$ by:
\begin{itemize}
\item $a_{1}=N-2$ and for $n\geq2$ by $a_{n}=a_{n-1}-\gcd(a_{n-1},N-(-1)^{n}(N-n))$ 
\end{itemize}
Then we claim:
\begin{itemize}
\item there is a least $g_{N}\in\left\{ 2,3,...,N-1\right\} $ such that
$a_{g_{N}}=0$ .
\item $g_{N}\sim N$ $(N\rightarrow\infty)$ 
\item for $N\geq2208$ we have $g_{N}+1$ and $2N-g_{N}-1$ which are simultanously
primes.$ $ 
\end{itemize}
Thus the Goldbach conjecture would be true. 

Hereafter a graph of $\frac{g_{N}}{N}\delta(g_{N}+1)\delta(2N-g_{N}-1)$
for $2\leq N\leq30000$. 

\begin{center}
(fig.14)
\par\end{center}

\begin{center}
\includegraphics[width=0.5\paperwidth]{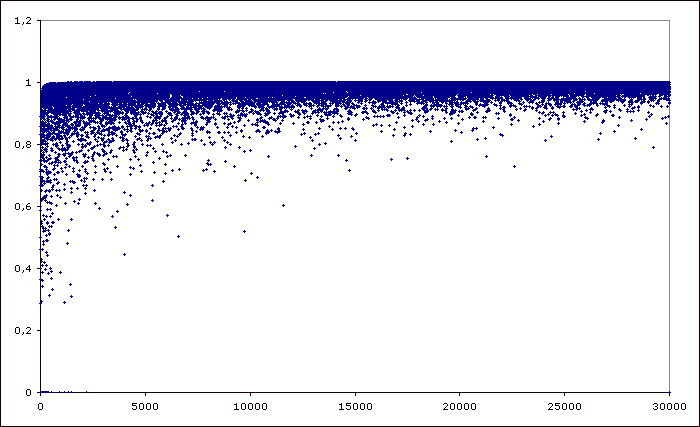}
\par\end{center}

We see there is no more zero value for $N\geq2208$. 

Next we provide a zoom of this picture in the range $100000\leq N\leq130000$. 

\begin{center}
(fig.15)
\par\end{center}

\begin{center}
\includegraphics[width=0.5\paperwidth]{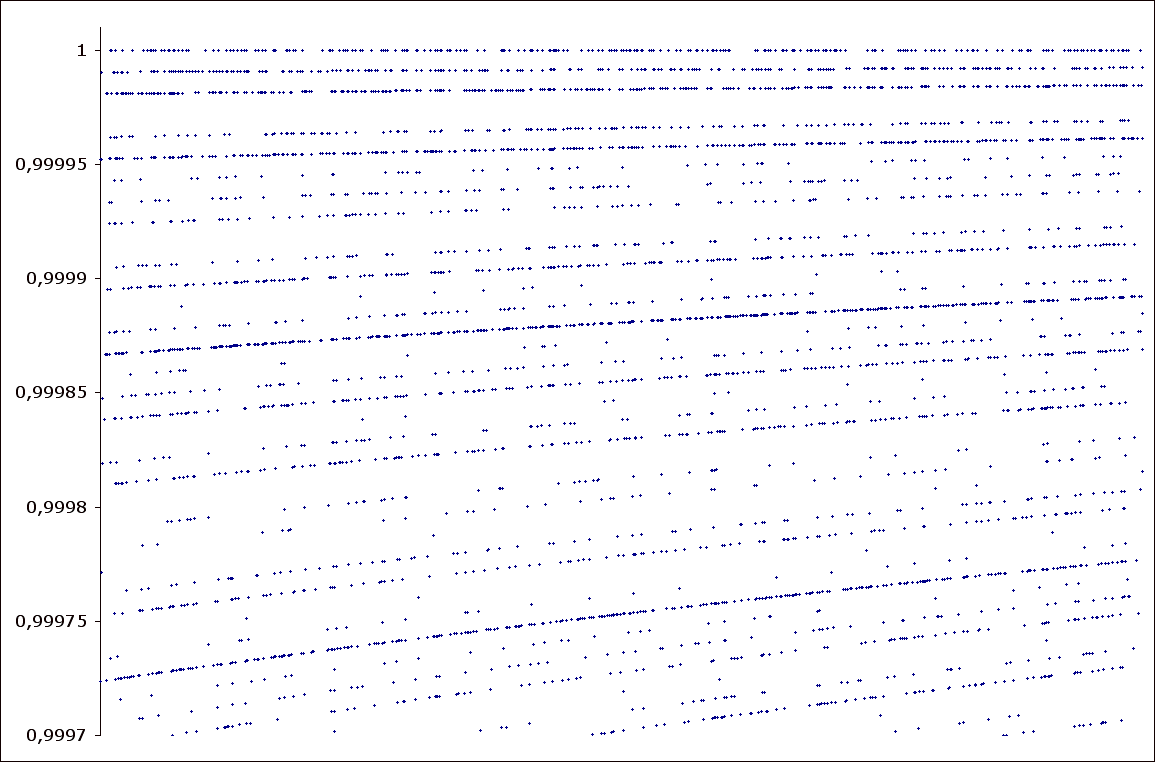}
\par\end{center}

The lines appearing regularly are due to the following observation.

\subsection*{Relation with the weak twin prime conjecture}

Moreover we claim that $\forall m\geq2\in\mathbb{N}$
\begin{itemize}
\item $g_{N}=N-m$ for infinitely many values of $N$. 
\end{itemize}
And we add withouth any other condition on $N$:
\begin{itemize}
\item $g_{N}=N-m\Leftrightarrow N-m+1$ and $N+m-1$ are both primes.
\end{itemize}
Thus for $m=2$ this means there are infinitely many twin primes.

\subsection*{The Goldbach constant}

It is worth to consider analytic aspects of this formulation of the
Goldbach conjecture. For instance we claim:
\begin{itemize}
\item $\sum_{k=3}^{n}\left(1-\frac{g_{k}}{k-1}\right)\sim C\sqrt{n}$ $\left(n\rightarrow\infty\right)$
where $C\leq4$. 
\end{itemize}
Hereafter the graph of $n^{-1/2}\sum_{k=3}^{n}\left(1-\frac{g_{k}}{k-1}\right)$
for $1\leq n\leq2000000$ 

\begin{center}
(fig.16)
\par\end{center}

\begin{center}
\includegraphics[width=0.5\paperwidth]{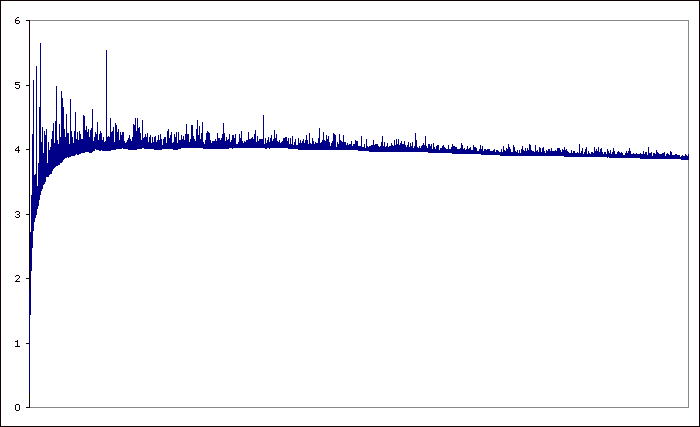}
\par\end{center}

\section{On a conjecture of Legendre}

In this study we give new formulation using gcd-algorithms for 3 problems
among a list of 4 problems considered by Landau in 1912 as \textquotedblleft{}unattackable
at the present state of science\textquotedblright{}\cite{key-6} \cite{key-7}:
\begin{itemize}
\item are they infinitely many primes numbers of form $n^{2}+1$?
\item the twin prime conjecture.
\item the Goldbach conjecture.
\end{itemize}
So a Landau problem is missing but our gcd-algorithm formulation still
works for this fourth probem attributed to Legendre. This remaining
conjecture says:
\begin{itemize}
\item there is always a prime between 2 consecutive squares. 
\end{itemize}
Our belief is that for $N$ large enough the algorithm given in 6.1
starting with $a(1)=(N+1)^{2}-2$ produces a prime value $f((N+1)^{2})$
greater than $N^{2}$. But this comes not from our asymptotic claim
in 6.1.

\subsubsection*{Remark}

Of course the situtation is clearer if we use the sequence $p(\pi(N))$
introduced in 6.1. and the conjectured asmptotic behaviour which works
for proving the Legendre conjecture (the behaviour toward zero of
$\frac{N-p(\pi(N))}{\sqrt{N}}$ shown in 6.1. ). But for this one
one must consider RH or things like that. Our hope is that the sequence
$f(N)$ despite its random nature allow us to avoid such considerations.

\subsection*{A Goldbach variation}

We came also across a stronger hypothesis than the original Legendre
conjecture and it is worth to present this {}``Goldbach variation''.
We don't know if this last conjecture could be helpful for finding
a {}``gcd-formulation'' of Legendre conjecture like our conjecture
10 for the Goldbach conjecture (we found nothing clear even for $N$
large enough). Namely we claim:
\begin{itemize}
\item $\forall N\geq175$ there is at least one $k$ satisfying $2\leq k\leq2N$
such that $N^{2}+k+1$ and $(N+1)^{2}-k$ are simultanously primes.
\end{itemize}
Here the graph of the function of $N$ for $2\leq N\leq20000$ representing
the number of $k$ verifying $2\leq k\leq2N$ and such that $N^{2}+k+1$
and $(N+1)^{2}-k$ are simultanously primes. 

\begin{center}
(fig.17)
\par\end{center}

\begin{center}
\includegraphics[width=0.5\paperwidth]{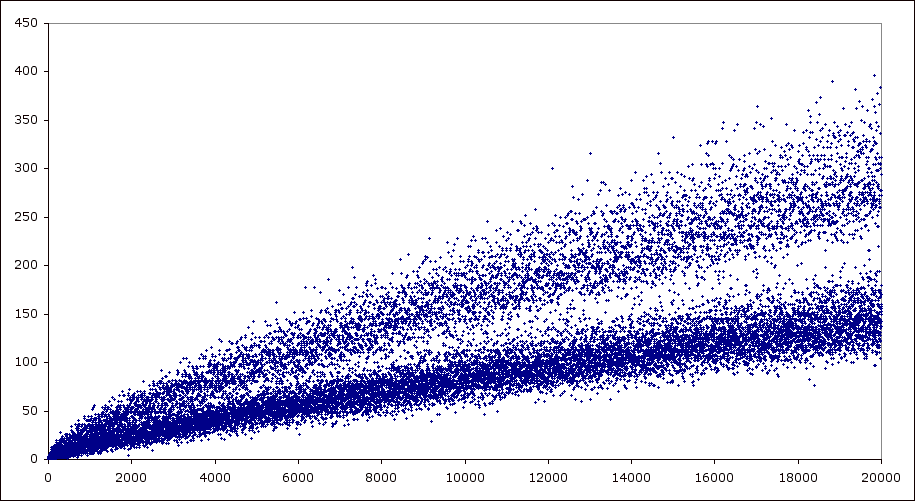}
\par\end{center}

It looks like the classical Goldbach comet.

\section*{Concluding remark}

The Schinzel hypothesis doesn't cover the Goldbach conjecture but
in view of our conjecture 10 we think there is a possible extension
of the conjecture 6. We mention that some people didn't think Rowland
simple theorem and related ideas of V. Shevelev could lead to interesting
results regarding the theory of the prime numbers. We hope we show
this advice is wrong and that some very important rules need to be
better understood and that our results could inspired some analytic
or probabilistic studies. Regarding the probabilistic approach we
think we have other striking arguments. In our study we consider only
natural objects in the gcd such as polynomials or periodic sequences.
It appears we can play with less predictive functions and for instance
we observe we can generate twin primes with any function $r$ taking
values in $\left\{ 0,2\right\} $ sufficiently {}``randomly''.%
\footnote{This needs to be clarified but we think $\lim_{n\rightarrow\infty}\frac{1}{n}\sum_{j=1}^{n}r_{j}\notin\left\{ 0,2\right\} $
is a sufficient condition and we believe the random process comes
mainly from the gcd.%
} To see this let us consider the differences of the Beatty sequence
for $\pi$ i.e. $v_{n}=\left\lfloor \pi n\right\rfloor -\left\lfloor \pi(n-1)\right\rfloor $
which takes values in $\left\{ 3,4\right\} $ and is not a periodic
sequence (this is the sequence A063438 in \cite{key-2}).

Then let $r_{n}=2(v_{n}-3)$ wich takes values in $\left\{ 0,2\right\} $
and define the sequence:
\begin{itemize}
\item $a(1)=N-2\geq0$
\item $a(n)=a(n-1)-\gcd(a(n-1),n+r_{n})$ 
\end{itemize}
Then we claim:
\begin{enumerate}
\item $\forall N\geq2$, $\exists f(N)\in\left\{ 1,2,\ldots,N\right\} $
such that $a(f(N))=0$.
\item $f(N)\sim N$ $\left(N\rightarrow\infty\right)$.
\item $\forall N\geq1649$, $f(N)+1$ and $f(N)+3$ are necessarily simultanously
primes.
\end{enumerate}
Hereafter we plot $\frac{f(N)}{N}\delta(f(N)+1)\delta(f(N)+3)$ for
$2\leq N\leq30000$. 

\begin{center}
(fig.18)
\par\end{center}

\begin{center}
\includegraphics[width=0.5\paperwidth]{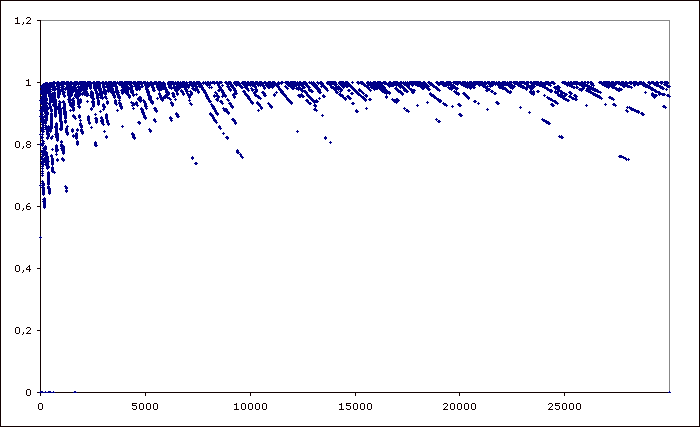}
\par\end{center}

We can see there is no more zero value for $N\geq1649$. 

\newpage{}

\newpage{}

\begin{center}
\textbf{\Large APPENDIX 1}
\par\end{center}{\Large \par}

\begin{center}
Table related to the conjecture 1
\par\end{center}

\begin{center}
\begin{tabular}{|c|c|c|c|}
\hline 
$k$ & $b_{1}(k)$ & $\delta(b_{1}(k))$ & $b_{1}(k)2^{-k}$\tabularnewline
\hline
\hline 
1 & 2  & 1 & 1.000000000\tabularnewline
\hline 
2 & 5  & 1 & 1.250000000\tabularnewline
\hline 
3 & 11 & 1 & 1.375000000\tabularnewline
\hline 
4 & 23 & 1 & 1.437500000\tabularnewline
\hline 
5 & 47 & 1 & 1.468750000\tabularnewline
\hline 
6 & 79 & 1 & 1.234375000\tabularnewline
\hline 
7 & 157 & 1 & 1.226562500\tabularnewline
\hline 
8 & 313 & 1 & 1.222656250\tabularnewline
\hline 
9 & 619 & 1 & 1.208984375\tabularnewline
\hline 
10 & 1237 & 1 & 1.208007812\tabularnewline
\hline 
11 & 2473 & 1 & 1.207519531\tabularnewline
\hline 
12 & 4909 & 1 & 1.198486328\tabularnewline
\hline 
13 & 9817 & 1 & 1.198364257\tabularnewline
\hline 
14 & 19603 & 1 & 1.196472167\tabularnewline
\hline 
15 & 39199 & 1 & 1.196258544\tabularnewline
\hline 
16 & 78193 & 1 & 1.193130493\tabularnewline
\hline 
17 & 156019 & 1 & 1.190330505\tabularnewline
\hline 
18 & 311347 & 1 & 1.187694549\tabularnewline
\hline 
19 & 622669 & 1 & 1.187646865\tabularnewline
\hline 
20 & 1244149 & 1 & 1.186512947\tabularnewline
\hline 
21 & 2487739 & 1 & 1.186246395\tabularnewline
\hline 
22 & 4975111 & 1 & 1.186158895\tabularnewline
\hline 
23 & 9950221 & 1 & 1.186158776\tabularnewline
\hline 
24 & 19900399 & 1 & 1.186156213\tabularnewline
\hline 
25 & 39800797 & 1 & 1.186156183\tabularnewline
\hline 
26 & 79601461 & 1 & 1.186154201\tabularnewline
\hline 
27 & 159202369 & 1 & 1.186150081\tabularnewline
\hline 
28 & 318404629 & 1 & 1.186149675\tabularnewline
\hline 
29 & 636788881 & 1 & 1.186111720\tabularnewline
\hline
\end{tabular}
\par\end{center}

\newpage{}

\begin{center}
\textbf{\Large APPENDIX 2}
\par\end{center}{\Large \par}

In the following table one can see the number of iterations $n$ giving
$a(n)=0$ (left column) and the corresponding value $mn+m-1$ in the
middle column. The right colum gives the value of the function {}``isprime''
using pari-gp.

\begin{center}
$m=5$
\par\end{center}

\begin{center}
\begin{tabular}{|c|c|c|}
\hline 
$n$ such that $a(n)=0$ & $5n+4$ & $\delta(5n+4)$\tabularnewline
\hline
\hline 
2 & 14 & 0\tabularnewline
\hline 
17 & 89 & 1\tabularnewline
\hline 
95 & 479 & 1\tabularnewline
\hline 
575 & 2879 & 1\tabularnewline
\hline 
3419 & 17099 & 1\tabularnewline
\hline 
19967 & 99839 & 1\tabularnewline
\hline 
119801 & 599009 & 1\tabularnewline
\hline 
718571 & 3592859 & 1\tabularnewline
\hline 
4311419 & 21557099 & 1\tabularnewline
\hline 
25867229 & 129336149 & 1\tabularnewline
\hline
\end{tabular}
\par\end{center}

\begin{center}
$m=6$
\par\end{center}

\begin{center}
\begin{tabular}{|c|c|c|}
\hline 
$n$ such that $a(n)=0$ & $6n+5$ & $\delta(6n+5)$\tabularnewline
\hline
\hline 
2 & 17 & 1\tabularnewline
\hline 
16 & 101 & 1\tabularnewline
\hline 
76 & 461 & 1\tabularnewline
\hline 
466 & 2801 & 1\tabularnewline
\hline 
3258 & 19553 & 1\tabularnewline
\hline 
22774 & 136649 & 1\tabularnewline
\hline 
159306 & 955841 & 1\tabularnewline
\hline 
1114124 & 6684749 & 1\tabularnewline
\hline 
77796204 & 46777229 & 1\tabularnewline
\hline 
54573434 & 327440609 & 1\tabularnewline
\hline
\end{tabular}
\par\end{center}

\begin{center}
$m=7$
\par\end{center}

\begin{center}
\begin{tabular}{|c|c|c|}
\hline 
$n$ such that $a(n)=0$ & $7n+6$ & $\delta(7n+6)$\tabularnewline
\hline
\hline 
2 & 21 & 0\tabularnewline
\hline 
23 & 167 & 1\tabularnewline
\hline 
113 & 797 & 1\tabularnewline
\hline 
899 & 6299 & 1\tabularnewline
\hline 
6973 & 48817 & 1\tabularnewline
\hline 
55633 & 389437 & 1\tabularnewline
\hline 
444901 & 3114313 & 1\tabularnewline
\hline 
3558575 & 24910031 & 1\tabularnewline
\hline 
28468585 & 199280101 & 1\tabularnewline
\hline
\end{tabular}
\par\end{center}

\begin{center}
$m=8$
\par\end{center}

\begin{center}
\begin{tabular}{|c|c|c|}
\hline 
$n$ such that $a(n)=0$ & $8n+7$ & $\delta(8n+7)$\tabularnewline
\hline
\hline 
2 & 23 & 1\tabularnewline
\hline 
20 & 167 & 1\tabularnewline
\hline 
188 & 1511 & 1\tabularnewline
\hline 
1682 & 13463 & 1\tabularnewline
\hline 
15020 & 120167 & 1\tabularnewline
\hline 
134504 & 1076039 & 1\tabularnewline
\hline 
1210544 & 9684359 & 1\tabularnewline
\hline 
10894874 & 87158999 & 1\tabularnewline
\hline 
98053784 & 784430279 & 1\tabularnewline
\hline
\end{tabular}
\par\end{center}

\begin{center}
$m=9$
\par\end{center}

\begin{center}
\begin{tabular}{|c|c|c|}
\hline 
$n$ such that $a(n)=0$ & $9n+8$ & $\delta(9n+8)$\tabularnewline
\hline
\hline 
2 & 26 & 0\tabularnewline
\hline 
29 & 269 & 1\tabularnewline
\hline 
299 & 2699 & 1\tabularnewline
\hline 
2935 & 26423 & 1\tabularnewline
\hline 
28869 & 259829 & 1\tabularnewline
\hline 
288385 & 2595473 & 1\tabularnewline
\hline 
2883809 & 25954289 & 1\tabularnewline
\hline 
28832339 & 259491059 & 1\tabularnewline
\hline
\end{tabular}
\par\end{center}

\begin{center}
$m=100$
\par\end{center}

\begin{center}
\begin{tabular}{|c|c|c|}
\hline 
$n$ such that $a(n)=0$ & $100n+99$ & $\delta(100n+99)$\tabularnewline
\hline
\hline 
2 & 299 & 0\tabularnewline
\hline 
226 & 22699 & 1\tabularnewline
\hline 
22810 & 2281099 & 1\tabularnewline
\hline 
2303908 & 230390899 & 1\tabularnewline
\hline
\end{tabular}
\par\end{center}

To see further how this is working nicely let us take $m=1000$. 

\begin{center}
\begin{tabular}{|c|c|c|}
\hline 
$n$ such that $a(n)=0$ & $1000n+999$ & $\delta(1000n+999)$\tabularnewline
\hline
\hline 
$2$ & $2999$ & $1$\tabularnewline
\hline 
$2986$ & $2986999$ & $1$\tabularnewline
\hline 
$2917174$ & $2917174999$ & $1$\tabularnewline
\hline
\end{tabular}
\par\end{center}

\paragraph{Finding big primes (from 1.3.1)}

We take $m=2^{k}$ for $100\leq k\leq200$ so that $b(1)=3$ and $a(b(1)+1)=3.2^{k}-1$
is the first record value (not necessarily prime of course and it
would be interesting to find conditions on the initial value in order
to have very few terms to compute). Next we stop the algorithm after
$\sim10^{5}$ iterations (few seconds are needed each time). In the
following table we keep the values of $k$ such that $(m+1)a(b(1)+1)+m+m\sum_{j=4}^{100000}\left(a(j+1)-a(j)+1\right)$
is a prime value of size $\sim3.4^{k}$and should be our second record
value $a(b(2)+1)$. We notice we came across 11 primes in that range. 

\begin{center}
\begin{tabular}{|c|c|}
\hline 
$k$ & $a(b(2)+1)$\tabularnewline
\hline
\hline 
100 & {\scriptsize 4820814132776970826625886270541990288599672051495735016816639}\tabularnewline
\hline 
107 & {\scriptsize 78984218751417890023438520762741628349070105828337814558023352319}\tabularnewline
\hline 
118 & {\scriptsize 331283824645947061796868281389238893531663193824530694464237235871940607}\tabularnewline
\hline 
127 & {\scriptsize 86844066927987146567678238756515930889442064948849015334398126094787194912767}\tabularnewline
\hline 
131 & {\scriptsize 22232081133564709521325629121668078276785918415807297904592970941163697569005567}\tabularnewline
\hline 
132 & {\scriptsize 88928324534258838085302516486672313230587227346404196132053363353462480809492479}\tabularnewline
\hline 
141 & {\scriptsize 23312026706708748851033542881882226879708773352492831418233933930240272261700160323583}\tabularnewline
\hline 
149 & {\scriptsize 1527776982250864564701334266307033620788600839243814368422923338893418735708806455350001663}\tabularnewline
\hline 
158 & {\scriptsize 400497569235170640449066569906791021488007090321237481912409760223102559592991367888359129087999}\tabularnewline
\hline 
164 & {\scriptsize 1640438043587258943279376670338216024014877072467663636588378426222267002018688816605037919597494271}\tabularnewline
\hline 
172 & {\scriptsize 107507747624534602106757229467285325349838983867263171688888770338883656854354650572386871166337445527551}\tabularnewline
\hline
\end{tabular}
\par\end{center}

\paragraph*{pari gp code}
\begin{description}
\item [{\texttt{\scriptsize for(k=100,200,m=2\textasciicircum{}k;a=1;S=0;M=0;}}]~{\scriptsize \par}
\item [{\texttt{\scriptsize for(n=2,10\textasciicircum{}5,t=a;a=abs(a-gcd(a,m{*}n-1));M=M+if(a-t>0,a-t,0);S=S+if(a-t<0,a-t+1,0);}}]~{\scriptsize \par}
\item [{\texttt{\scriptsize if(abs(t-a)>1,if(isprime((m+1){*}M+m+m{*}S)==1,print(k,''}}] \texttt{\scriptsize {}``,(m+1){*}M+m+m{*}S,\textquotedbl{}\textquotedbl{})))))}{\scriptsize \par}
\end{description}
\newpage{}

\begin{center}
\textbf{\Large APPENDIX 3}
\par\end{center}{\Large \par}

Let $a(1)=m$ and 
\begin{itemize}
\item $a(n)=\left|a(n-1)-\gcd(a(n-1),n^{2}-1)\right|.$ 
\end{itemize}
Here a table for various starting values $m$ allowing us to exhibit
distinct pairs of twin primes (the other values $m=10k$ for $1\leq k\leq14$
produce also pairs of twin primes but they are all in this list keeping
the distinct pairs only).

\begin{center}
\begin{tabular}{|c|c|c|c|}
\hline 
$m$ & $n$ such that $a(n)=0$ & $\delta(n)$  & $\delta(n+2)$\tabularnewline
\hline
\hline 
10  & 11  & 1 & 1\tabularnewline
\hline 
 & 137 & 1 & 1\tabularnewline
\hline 
 & 19181 & 1 & 1\tabularnewline
\hline 
20 & 17 & 1 & 1\tabularnewline
\hline 
 & 281 & 1 & 1\tabularnewline
\hline 
 & 79559 & 1 & 1\tabularnewline
\hline 
30 & 29 & 1 & 1\tabularnewline
\hline 
 & 881 & 1 & 1\tabularnewline
\hline 
 & 777011 & 1 & 1\tabularnewline
\hline 
40 & 41 & 1 & 1\tabularnewline
\hline 
 & 1787 & 1 & 1\tabularnewline
\hline 
 & 3198731 & 1 & 1\tabularnewline
\hline 
60 & 59 & 1 & 1\tabularnewline
\hline 
 & 3527 & 1 & 1\tabularnewline
\hline 
 & 12448001 & 1 & 1\tabularnewline
\hline 
70 & 71 & 1 & 1\tabularnewline
\hline 
 & 5099 & 1 & 1\tabularnewline
\hline 
 & 26010041 & 1 & 1\tabularnewline
\hline 
100 & 101 & 1 & 1\tabularnewline
\hline 
 & 10499 & 1 & 1\tabularnewline
\hline 
 & 110258891 & 1 & 1\tabularnewline
\hline 
110 & 107 & 1 & 1\tabularnewline
\hline 
 & 11699 & 1 & 1\tabularnewline
\hline 
 & 136890881 & 1 & 1\tabularnewline
\hline 
140 & 137 & 1 & 1\tabularnewline
\hline 
 & 19181 & 1 & 1\tabularnewline
\hline 
 & 367953497 & 1 & 1\tabularnewline
\hline
\end{tabular}
\par\end{center}

And 2 more starting values

\begin{center}
\begin{tabular}{|c|c|c|c|}
\hline 
$m$ & $n$ such that $a(n)=0$ & $\delta(n)$  & $\delta(n+2)$\tabularnewline
\hline
\hline 
200  & 197 & 1 & 1\tabularnewline
\hline 
 & 39161 & 1 & 1\tabularnewline
\hline 
 & 1533646397 & 1 & 1\tabularnewline
\hline 
300 & 281 & 1 & 1\tabularnewline
\hline 
 & 79559 & 1 & 1\tabularnewline
\hline 
 & 6329815697 & 1 & 1\tabularnewline
\hline
\end{tabular}
\par\end{center}

We also provide this computation but using Rowland-Shevelev recursion:
\begin{itemize}
\item $a_{1}=3$ then $a_{n}=a_{n-1}+\gcd(a_{n-1,}n(n-2))$ 
\end{itemize}
Here the table of records of the differences $a_{n}-a_{n-1}$ which
are the lower terms of some twin prime pairs.

\begin{center}
\begin{tabular}{|c|}
\hline 
records of $a_{n}-a_{n-1}$\tabularnewline
\hline
\hline 
3\tabularnewline
\hline 
5\tabularnewline
\hline 
11\tabularnewline
\hline 
41\tabularnewline
\hline 
101\tabularnewline
\hline 
239\tabularnewline
\hline 
521\tabularnewline
\hline 
1049\tabularnewline
\hline 
2111\tabularnewline
\hline 
4229\tabularnewline
\hline 
10331 \tabularnewline
\hline 
20747 \tabularnewline
\hline 
41519 \tabularnewline
\hline 
83219 \tabularnewline
\hline 
166847 \tabularnewline
\hline 
333791 \tabularnewline
\hline 
669479 \tabularnewline
\hline 
1341017 \tabularnewline
\hline 
2682539 \tabularnewline
\hline 
5365229 \tabularnewline
\hline 
10732751 \tabularnewline
\hline 
21466259 \tabularnewline
\hline 
42932567 \tabularnewline
\hline 
85865321 \tabularnewline
\hline 
171730679 \tabularnewline
\hline 
343461647 \tabularnewline
\hline 
686929511 \tabularnewline
\hline 
1373891861 \tabularnewline
\hline 
2747784329 \tabularnewline
\hline 
5495586839 \tabularnewline
\hline
\end{tabular}
\par\end{center}

\newpage{}

\begin{center}
\textbf{\Large APPENDIX 4}
\par\end{center}{\Large \par}

Here $a(1)=4m^{2}$ and $a(n)=\left|a(n-1)-\gcd(a(n-1),n(n+2m))\right|$
and the 3 first values of $n$ such that $a_{n}=0$ .

\begin{center}
\begin{tabular}{|c|c|c|c|}
\hline 
$m$ & $n+1$ when $a(n)=0$ & $\delta(n+1)$  & $\delta(n+2m+1)$\tabularnewline
\hline
\hline 
2  & 13 & 1 & 1\tabularnewline
\hline 
 & 193 & 1 & 1\tabularnewline
\hline 
 & 38197 & 1 & 1\tabularnewline
\hline 
3 & 31 & 1 & 1\tabularnewline
\hline 
 & 1091 & 1 & 1\tabularnewline
\hline 
 & 1197193 & 1 & 1\tabularnewline
\hline 
4 & 59 & 1 & 1\tabularnewline
\hline 
 & 4013 & 1 & 1\tabularnewline
\hline 
 & 16138511 & 1 & 1\tabularnewline
\hline 
5 & 97 & 1 & 1\tabularnewline
\hline 
 & 10477 & 1 & 1\tabularnewline
\hline 
 & 109880317 & 1 & 1\tabularnewline
\hline 
6 & 139 & 1 & 1\tabularnewline
\hline 
 & 20521 & 1 & 1\tabularnewline
\hline 
 & 421370778 & 1 & 1\tabularnewline
\hline
\end{tabular}
\par\end{center}

\newpage{}

\begin{center}
\textbf{\Large APPENDIX 5}
\par\end{center}{\Large \par}

\begin{center}
$p=2$
\par\end{center}

\begin{center}
\begin{tabular}{|c|c|c|}
\hline 
Values of $n$ such that $a(n)=0$ & $2(n+1)^{2}-1$ & $\delta(2(n+1)^{2}-1)$\tabularnewline
\hline
\hline 
3 & 31 & 1\tabularnewline
\hline 
35 & 2591 & 1\tabularnewline
\hline 
2627 & 13812767 & 1\tabularnewline
\hline 
11993333 & 287680120871111 & 1\tabularnewline
\hline
\end{tabular}
\par\end{center}

\begin{center}
$p=3$
\par\end{center}

\begin{center}
\begin{tabular}{|c|c|c|}
\hline 
Values of $n$ such that $a(n)=0$ & $3(n+1)^{2}-1$ & $\delta(3(n+1)^{2}-1)$\tabularnewline
\hline
\hline 
3 & 47 & 1\tabularnewline
\hline 
51 & 8111 & 1\tabularnewline
\hline 
7665 & 176302667 & 1\tabularnewline
\hline 
176310323 & 93255991046954927 & 1\tabularnewline
\hline
\end{tabular}
\par\end{center}

\begin{center}
$p=5$
\par\end{center}

\begin{center}
\begin{tabular}{|c|c|c|}
\hline 
Values of $n$ such that $a(n)=0$ & $5(n+1)^{2}-1$ & $\delta(5(n+1)^{2}-1)$\tabularnewline
\hline
\hline 
3 & 79 & 1\tabularnewline
\hline 
83 & 35279 & 1\tabularnewline
\hline 
34647 & 6002419519 & 1\tabularnewline
\hline
\end{tabular}
\par\end{center}

\begin{center}
$p=7$
\par\end{center}

\begin{center}
\begin{tabular}{|c|c|c|}
\hline 
Values of $n$ such that $a(n)=0$ & $7(n+1)^{2}-1$ & $\delta(7(n+1)^{2}-1)$\tabularnewline
\hline
\hline 
5 & 251 & 1\tabularnewline
\hline 
257 & 465947 & 1\tabularnewline
\hline 
461009 & 1487711540699 & 1\tabularnewline
\hline
\end{tabular}
\par\end{center}

\begin{center}
$p=11$
\par\end{center}

\begin{center}
\begin{tabular}{|c|c|c|}
\hline 
Values of $n$ such that $a(n)=0$ & $11(n+1)^{2}-1$ & $\delta(11(n+1)^{2}-1)$\tabularnewline
\hline
\hline 
11 & 1583 & 1\tabularnewline
\hline 
1419 & 22180399 & 1\tabularnewline
\hline 
22181509 & 5412213244681099 & 1\tabularnewline
\hline
\end{tabular}
\par\end{center}

\begin{center}
$p=17$
\par\end{center}

\begin{center}
\begin{tabular}{|c|c|c|}
\hline 
Values of $n$ such that $a(n)=0$ & $17(n+1)^{2}-1$ & $\delta(17(n+1)^{2}-1)$\tabularnewline
\hline
\hline 
11 & 2447 & 1\tabularnewline
\hline 
2417 & 99394307 & 1\tabularnewline
\hline 
87543523 & 130285766103755791 & 1\tabularnewline
\hline
\end{tabular}
\par\end{center}

\newpage{}

\begin{center}
\textbf{\Large APPENDIX 6}
\par\end{center}{\Large \par}

We define the sequence $a$ as follows:
\begin{itemize}
\item $a(1)=k$ 
\item $a(n)=\left|a(n-1)-\gcd(a(n-1),mn+b(n))\right|$ 
\end{itemize}
We then give tables supporting the conjecture 9 for various $k,m$
and $b$ periodic sequence which is given by its period $\left\{ b_{1},b_{2},...,b_{\beta}\right\} $.

\paragraph*{2 periodic sequence}

\begin{center}
$b=\left\{ 0,2\right\} ,m=1,k=100$
\par\end{center}

\begin{center}
\begin{tabular}{|c|c|}
\hline 
$n$ such that $a(n)=0$ & $\delta(n+1)\delta(n+3)$\tabularnewline
\hline
\hline 
100 & 1\tabularnewline
\hline 
196 & 1\tabularnewline
\hline 
310 & 1\tabularnewline
\hline 
616 & 1\tabularnewline
\hline 
1228 & 1\tabularnewline
\hline 
2380 & 1\tabularnewline
\hline 
4648 & 1\tabularnewline
\hline 
8860 & 1\tabularnewline
\hline 
17026 & 1\tabularnewline
\hline 
33808  & 1\tabularnewline
\hline 
67408  & 1\tabularnewline
\hline 
134680 & 1\tabularnewline
\hline 
267718 & 1\tabularnewline
\hline 
535348  & 1\tabularnewline
\hline 
1069216 & 1\tabularnewline
\hline 
2138398  & 1\tabularnewline
\hline 
4275640  & 1\tabularnewline
\hline 
8545696  & 1\tabularnewline
\hline 
17091376  & 1\tabularnewline
\hline 
34182748  & 1\tabularnewline
\hline 
68365468  & 1\tabularnewline
\hline 
136730638  & 1\tabularnewline
\hline 
273461158  & 1\tabularnewline
\hline 
546917140  & 1\tabularnewline
\hline 
1093813726  & 1\tabularnewline
\hline 
2187610990 & 1\tabularnewline
\hline
\end{tabular}
\par\end{center}

\begin{center}
\newpage{}$b=\left\{ -1,1\right\} ,m=10,k=100$
\par\end{center}

\begin{center}
\begin{tabular}{|c|c|}
\hline 
$n$ such that $a(n)=0$ & $\delta(10n+9)\delta(10n+11)$\tabularnewline
\hline
\hline 
101 & 1\tabularnewline
\hline 
1115 & 1\tabularnewline
\hline 
12203 & 1\tabularnewline
\hline 
130013 & 1\tabularnewline
\hline 
1427183  & 1\tabularnewline
\hline 
15692309 & 1\tabularnewline
\hline 
172614683 & 1\tabularnewline
\hline
\end{tabular}
\par\end{center}

\begin{center}
$b=\left\{ -2,2\right\} ,m=3,k=100$
\par\end{center}

\begin{center}
\begin{tabular}{|c|c|}
\hline 
$n$ such that $a(n)=0$ & $\delta(3n+1)\delta(3n+5)$\tabularnewline
\hline
\hline 
34 & 1\tabularnewline
\hline 
116 & 1\tabularnewline
\hline 
434 & 1\tabularnewline
\hline 
1576 & 1\tabularnewline
\hline 
6102 & 1\tabularnewline
\hline 
21154  & 1\tabularnewline
\hline 
84606  & 1\tabularnewline
\hline 
338386  & 1\tabularnewline
\hline 
1351382  & 1\tabularnewline
\hline 
5405526  & 1\tabularnewline
\hline 
21622094  & 1\tabularnewline
\hline
\end{tabular}
\par\end{center}

\paragraph*{3 periodic sequence}

\begin{center}
$b=\left\{ 2,6,0\right\} ,m=1,k=3000$
\par\end{center}

\begin{center}
\begin{tabular}{|c|c|}
\hline 
$n$ such that $a(n)=0$ & $\delta(n+1)\delta(n+3)\delta(n+7)$\tabularnewline
\hline
\hline 
2686 & 1\tabularnewline
\hline 
5230 & 1\tabularnewline
\hline 
10456  & 1\tabularnewline
\hline 
19420 & 1\tabularnewline
\hline 
29566 & 1\tabularnewline
\hline 
54496  & 1\tabularnewline
\hline 
105526  & 1\tabularnewline
\hline 
211060  & 1\tabularnewline
\hline 
408430  & 1\tabularnewline
\hline 
802126  & 1\tabularnewline
\hline 
1600216  & 1\tabularnewline
\hline 
3200200  & 1\tabularnewline
\hline 
6393910 & 1\tabularnewline
\hline 
12783496 & 1\tabularnewline
\hline 
25566676 & 1\tabularnewline
\hline 
51095410 & 1\tabularnewline
\hline 
102190390 & 1\tabularnewline
\hline 
204347176 & 1\tabularnewline
\hline
\end{tabular}
\par\end{center}

\begin{center}
$b=\left\{ 0,4,6\right\} ,m=1,k=3000$
\par\end{center}

\begin{center}
\begin{tabular}{|c|c|}
\hline 
$n$ such that $a(n)=0$ & $\delta(n+1)\delta(n+5)\delta(n+7)$\tabularnewline
\hline
\hline 
2082 & 1\tabularnewline
\hline 
3462 & 1\tabularnewline
\hline 
6546 & 1\tabularnewline
\hline 
12372  & 1\tabularnewline
\hline 
23052  & 1\tabularnewline
\hline 
44262  & 1\tabularnewline
\hline 
85086  & 1\tabularnewline
\hline 
167016  & 1\tabularnewline
\hline 
313986  & 1\tabularnewline
\hline 
622476  & 1\tabularnewline
\hline 
1237206  & 1\tabularnewline
\hline 
2452752  & 1\tabularnewline
\hline 
4882326  & 1\tabularnewline
\hline 
9753276  & 1\tabularnewline
\hline 
19504866  & 1\tabularnewline
\hline
\end{tabular}
\par\end{center}

\begin{center}
$b=\left\{ -2,2,-4\right\} ,m=5,k=2000$
\par\end{center}

\begin{center}
\begin{tabular}{|c|c|}
\hline 
$n$ such that $a(n)=0$ & $\delta(5n+1)\delta(5n+3)\delta(5n+7)$\tabularnewline
\hline
\hline 
1772 & 1\tabularnewline
\hline 
9806 & 1\tabularnewline
\hline 
58274 & 1\tabularnewline
\hline 
343772  & 1\tabularnewline
\hline 
2057378  & 1\tabularnewline
\hline 
12342518  & 1\tabularnewline
\hline 
73895180 & 1\tabularnewline
\hline
\end{tabular}
\par\end{center}

\newpage{}

\paragraph*{4 periodic sequence}

\begin{center}
$b=\left\{ 1,7,11,-1\right\} ,m=1,k=20000$
\par\end{center}

\begin{center}
\begin{tabular}{|c|c|}
\hline 
$n$ such that $a(n)=0$ & $\delta(n)\delta(n+2)\delta(n+8)\delta(n+12)$\tabularnewline
\hline
\hline 
19421 & 1\tabularnewline
\hline 
36779 & 1\tabularnewline
\hline 
70841  & 1\tabularnewline
\hline 
 138239 & 1\tabularnewline
\hline 
 236771 & 1\tabularnewline
\hline 
 443159 & 1\tabularnewline
\hline 
 882239 & 1\tabularnewline
\hline 
 1758389 & 1\tabularnewline
\hline 
 3376979 & 1\tabularnewline
\hline 
 6631901 & 1\tabularnewline
\hline 
13236539  & 1\tabularnewline
\hline 
 26425379  & 1\tabularnewline
\hline 
 52658999 & 1\tabularnewline
\hline 
 104785649 & 1\tabularnewline
\hline 
209560319 & 1\tabularnewline
\hline 
418973999 & 1\tabularnewline
\hline
\end{tabular}
\par\end{center}

\begin{center}
$b=\left\{ 2,6,18,26\right\} ,m=1,k=100000$
\par\end{center}

\begin{center}
\begin{tabular}{|c|c|}
\hline 
$n$ such that $a(n)=0$ & $\prod_{j=1}^{4}\delta(n+b(j)+1)$\tabularnewline
\hline
\hline 
83200 & 1\tabularnewline
\hline 
150190 & 1\tabularnewline
\hline 
294754  & 1\tabularnewline
\hline 
573844  & 1\tabularnewline
\hline 
1107784  & 1\tabularnewline
\hline 
2208064 & 1\tabularnewline
\hline 
4171774  & 1\tabularnewline
\hline 
8332840 & 1\tabularnewline
\hline 
16461094 & 1\tabularnewline
\hline 
32756680  & 1\tabularnewline
\hline 
65166814  & 1\tabularnewline
\hline 
130175344  & 1\tabularnewline
\hline 
260331034  & 1\tabularnewline
\hline 
520484380  & 1\tabularnewline
\hline 
1040389234  & 1\tabularnewline
\hline 
2080515244  & 1\tabularnewline
\hline 
4161006904 & 1\tabularnewline
\hline 
8321226490 & 1\tabularnewline
\hline
\end{tabular}
\par\end{center}

\begin{center}
\newpage{}
\par\end{center}

\paragraph*{5 periodic sequence}

\begin{center}
$b=\left\{ 2,8,12,14,18\right\} ,m=1,k=1500000$
\par\end{center}

\begin{center}
\begin{tabular}{|c|c|}
\hline 
$n$ such that $a(n)=0$ & $\prod_{j=1}^{5}\delta(n+1+b(j))$\tabularnewline
\hline
\hline 
1212424 & 1\tabularnewline
\hline 
2270674 & 1\tabularnewline
\hline 
4271158 & 1\tabularnewline
\hline 
8358658 & 1\tabularnewline
\hline 
15875398 & 1\tabularnewline
\hline 
31562608 & 1\tabularnewline
\hline 
62555878 & 1\tabularnewline
\hline 
125087098 & 1\tabularnewline
\hline 
249509788 & 1\tabularnewline
\hline 
477331018 & 1\tabularnewline
\hline 
954642034 & 1\tabularnewline
\hline 
1905881278 & 1\tabularnewline
\hline 
3809937208 & 1\tabularnewline
\hline
\end{tabular}
\par\end{center}

\begin{center}
$b=\left\{ 2,8,12,14,18\right\} ,m=1,k=2000000$
\par\end{center}

\begin{center}
\begin{tabular}{|c|c|}
\hline 
$n$ such that $a(n)=0$ & $\prod_{j=1}^{5}\delta(n+1+b(j))$\tabularnewline
\hline
\hline 
1460728 & 1\tabularnewline
\hline 
2839924 & 1\tabularnewline
\hline 
4218154 & 1\tabularnewline
\hline 
8068438  & 1\tabularnewline
\hline 
16130884 & 1\tabularnewline
\hline 
32240278  & 1\tabularnewline
\hline 
64123234  & 1\tabularnewline
\hline 
127725328 & 1\tabularnewline
\hline 
254416288  & 1\tabularnewline
\hline 
507764278 & 1\tabularnewline
\hline
\end{tabular}
\par\end{center}

\paragraph*{6 periodic sequence}

\begin{center}
$b=\left\{ 0,2,6,14,30,62\right\} ,m=1,k=2000000$
\par\end{center}

\begin{center}
\begin{tabular}{|c|c|}
\hline 
$n$ such that $a(n)=0$ & $\prod_{j=1}^{6}\delta(n+1+b(j))$\tabularnewline
\hline
\hline 
1460728 & 1\tabularnewline
\hline 
2839924 & 1\tabularnewline
\hline 
4218154 & 1\tabularnewline
\hline 
8068438  & 1\tabularnewline
\hline 
16130884 & 1\tabularnewline
\hline 
32240278  & 1\tabularnewline
\hline 
64123234  & 1\tabularnewline
\hline 
127725328 & 1\tabularnewline
\hline 
254416288  & 1\tabularnewline
\hline 
507764278 & 1\tabularnewline
\hline
\end{tabular}
\par\end{center}
\end{document}